\documentclass[12pt]{amsart}

\numberwithin{equation}{section}
\newtheorem{theorem}{Theorem}[section]
\newtheorem{lemma}[theorem]{Lemma}
\newtheorem{proposition}[theorem]{Proposition}
\newtheorem{corollary}[theorem]{Corollary}
\newtheorem*{corollary*}{Corollary}
\newtheorem*{proposition*}{Proposition}
\newtheorem*{theorem*}{Theorem}
\newtheorem*{question}{Question}
\theoremstyle{definition}
\newtheorem*{definition}{Definition}
\theoremstyle{plain}
\newtheorem*{Szth}{Szemer\'edi's Theorem}
\newtheorem*{Fth}{Furstenberg's Theorem}
\newtheorem*{CorrPrinc}{Furstenberg's Correspondence Principle}
\newtheorem*{PMRT}{Polynomial Multiple Recurrence Theorem}
\newtheorem*{PST}{Polynomial Szemer\'edi's Theorem}
\newtheorem*{MRT}{Multidimensional Recurrence Theorem}
\newtheorem*{MST}{Multidimensional Szemer\'edi's Theorem}
\newtheorem*{VdC}{Van der Corput Lemma}
\newtheorem*{ST}{Structure Theorem}

\newcommand{\C}{{\mathbb C}}
\newcommand{\E}{{\mathbb E}}
\newcommand{\N}{{\mathbb N}}
\newcommand{\Q}{{\mathbb Q}}
\newcommand{\R}{{\mathbb R}}
\newcommand{\T}{{\mathbb T}}
\newcommand{\Z}{{\mathbb Z}}
\newcommand{\CC}{{\mathcal C}}
\newcommand{\CG}{{\mathcal G}}
\newcommand{\CH}{{\mathcal H}}
\newcommand{\CI}{{\mathcal I}}
\newcommand{\CJ}{{\mathcal J}}
\newcommand{\CS}{{\mathcal S}}
\newcommand{\CW}{{\mathcal W}}
\newcommand{\CX}{{\mathcal X}}
\newcommand{\CY}{{\mathcal Y}}
\newcommand{\CZ}{{\mathcal Z}}
\newcommand{\ba}{\mathbf{a}}
\newcommand{\bn}{\mathbf{n}}
\newcommand{\bt}{\mathbf{t}}
\newcommand{\bu}{\mathbf{u}}
\newcommand{\bx}{\mathbf{x}}
\newcommand{\by}{\mathbf{y}}
\newcommand{\zero}{\mathbf{0}}
\newcommand{\one}{\mathbf{1}}
\newcommand{\norm}[1]{\lVert #1\rVert}
\newcommand{\nnorm}[1]{\lvert\!|\!| #1|\!|\!\rvert}
\DeclareMathOperator{\card}{Card}
\DeclareMathOperator{\id}{Id}
\newcommand{\inv}{^{-1}}
\newcommand{\eps}{\varepsilon}
\newcommand{\type}[1]{^{[#1]}}
\newcommand{\typ}[1]{^{[#1]^*}}
\begin{document}

\title{Convergence of multiple ergodic averages}
\author{Bernard Host}
\address{Math\'ematiques, Universit\'e de Marne la Vall\'ee.
5, Bd. Descartes, Champs sur Marne.
77454 Marne la Vall\'ee Cedex, France}
\email{host@math.univ-mlv.fr}

\begin{abstract}
These notes are based on a course for a general audience 
given at the Centro de Modeliamento 
Matem\'atico of the University of Chile, in December 2004. 

We study the mean convergence of multiple ergodic averages, that
is, averages of a product of functions taken at different times.
We also describe the relations between this area of ergodic theory
and some classical and some recent results in additive number 
theory.
\end{abstract}

\maketitle

In this paper we present some recent theorems of convergence of
multiple ergodic averages.
While the classical Ergodic Theorem describes the limit behavior of
the time averages of a function, these theorems deal with averages of
a product of functions taken at different times. We essentially focus on
the case that the $k$ functions are taken at times $n$, $2n$, \dots,
$kn$ (Theorem~\ref{th:convprog}) but we also consider the case of
polynomial times.

These convergence results belong to the field initiated
by Furstenberg,
exploring recurrence properties in ergodic theory and their relations
with combinatorial number theory.
The Correspondence Principle (see Section~\ref{sec:context}) provides
a bridge between these two domains by allowing one to deduce 
combinatorial properties of sets of integers with positive upper density
from recurrence theorems.
For example, Szemer\'edi's Theorem on the existence of arithmetic
progressions in sets of positive upper density corresponds to Furstenberg's
Theorem about multiple recurrence along arithmetic progressions.
The original proof of Szemer\'edi's Theorem is purely combinatorial
and some of its generalizations also
have  combinatorial proofs, apparently completely
different of the ergodic ones.
But recent progress in both fields leads to the intuition that there
exists a hidden relation between the objects and methods of the
two areas.
Understanding this relation more precisely could be an interesting
challenge. For this reason we do not completely ignore the
combinatorial point of  view in this paper.

As  often happens, the convergence theorems of multiple ergodic
averages will probably receive short self-contained proofs sometime in the
future.
But this is not the case  at the present  time
and it would be completely impossible to give complete proofs within the
framework of this paper.
Our more limited goal is to present most of the necessary tools
and to summarize the main steps.
Some partial results are given complete or partial proofs.
These are neither the most difficult nor the ``most important'',
but they are chosen for the enlightenment they  bring to some of the
main ingredients and to their uses.

We hope that the majority of this text is
accessible to a reader with a minimal knowledge in ergodic theory.
For  easier reading some notes have been postponed to the end of
the sections and further comments are included in the Appendix.
This material is not used in the main text and is intended as a
supplemental material.
The paper uses a lot of notation and some of it is not classical.
We tried to keep notation as similar as possible to the original
papers referred to.
For an easier reading, we introduce the notation
only when it is needed.
However it should be noted that throughout, we use the symbol $\N$
  with its European meaning:
$\N=\{0,1,2,\dots\}$.

\section{Context and  Results}
\label{sec:context}
\subsection{Szemer\'edi's Theorem}
Our starting point is a celebrated theorem by
Szemer\'edi. We begin with some definitions.
\begin{definition}
The \emph{upper density} of a subset $E$ of $\N$ is:
$$
    d^*(E)=\limsup_{N\to\infty}\frac 1N\card(E\cap\{0,1,2,\dots,N-1\})\ .
$$
An arithmetic progression of length $k$ is a set of integers of the
form
$$
\{a,a+d,\dots, a+(k-1)d\}
$$
    where $a,d$ are integers and $d>0$.
\end{definition}

\begin{Szth}[\cite{Sz2}]
\label{th:SZ}
A subset of integers with positive upper
density contains arbitrarily long  arithmetic progressions.
\end{Szth}

The result was first conjectured by Erd\"os and Turan~\cite{ET} in 1936
and was solved by Roth~\cite{Roth} in 1953 for progressions of length
$3$ and by Szemer\'edi~\cite{Sz1} in 1969 for progressions of length $4$.
While Roth's proof  belongs to harmonic analysis, Szemer\'edi's method
is combinatorial and relies on graph theory.

This theorem can be reformulated in terms of finite sets:

\begin{theorem}
\label{th:SZF}
For every integer $k\geq 2$ and every real $\delta>0$ there exists
an integer $N(k,\delta)$ such that:\\
For every $N>N(k,\delta)$, every subset $E$ of $\{1,2,\dots,N\}$
with $\card(E)>\delta N$ contains an arithmetic progression of length
$k$.
\end{theorem}

Szemer\'edi's Theorem follows immediately from
its finite version and the converse implication uses a simple
compactness argument. It is worth noting that this method cannot
provide an explicit value for the constants $N(k,\delta)$. The
original proof of Szemer\'edi did not give (usable) constants either
and many in the combinatorial community competed to
find the best constants for progressions of length $3$ and $4$ (the
winner for length $3$ was Bourgain~\cite{Bourgain2}) until a few years
ago when Gowers~\cite{Gowers1} gave a new proof with explicit
constants for the general case.  He proved the theorem in a more
analytical form (see~\ref{subsec:anal}),
already used by Roth and Bourgain:

\begin{theorem*}
Let $N\geq 2$ be an integer and let $\Z/N\Z$ be endowed with its
normalized Haar measure $m$. For
every integer $k\geq 2$ and every real $\delta>0$ there exists
a constant $c(k,\delta)>0$, not depending on $N$, such that:\\
For every function $f$ on $\Z/N\Z$ with $0\leq f\leq 1$ and $\int
f(x)\,dm(x)\geq\delta$,
\begin{equation}
\label{eq:SzAnal}
    \iint f(x)f(x+y)f(x+2y)\dots f(x+(k-1)y)\,dm(x)\,dm(y)\geq
    c(k,\delta)\ .
\end{equation}
\end{theorem*}

Gowers' proof uses methods of both Fourier Analysis (the circle
method) and combinatorial number theory, in particular Freiman's
results~(\cite{Fr1}, \cite{Fr2}) on the sums
of sets of integers.

Although  it does not seem related to our topic it would be difficult
not to mention the spectacular new result by Green and Tao which
answered a very old question.
\begin{theorem*}[Green \& Tao~\cite{GT}]
The  set of primes contains arbitrarily long arithmetic
progressions.
\end{theorem*}
(See also  Note~\ref{subsec:noteGT}.)
There exists a relation between the proof of Green and Tao and ergodic
theory, even though this relation is not completely understood at
this time (see
Appendix~\ref{subsec:gowersnorms}).

\subsection{Furstenberg's Theorem and its generalizations}
    Before stating the result we fix some notation and
some conventions.

In general, we write $(X,\mu)$ for a probability space, omitting the
$\sigma$-algebra;  when needed it is denoted by  by the corresponding
calligraphic letter $\CX$.
We always assume that  $\CX$ is countably generated and that
$\CX$ is the Borel $\sigma$-algebra whenever $X$ is
endowed with a (Polish) topology.
Throughout these notes, every subset of $X$ is implicitly assumed to be
    measurable and the term ``bounded
function on $X$'' means bounded and measurable.

By a \emph{system}, we mean a probability space $(X,\mu)$ endowed with
an invertible, bi-measurable, measure preserving transformation
$T\colon X\to X$ and we write the system as $(X, \mu, T)$. In the main
results the
hypothesis of invertibility of $T$ can be removed by passing to the
natural extension.

In 1977 Furstenberg  proved a very beautiful result about
multiple recurrence in ergodic theory:

\begin{Fth}[\cite{F2}; \cite{FKO} is easier to read]
Let  $(X,\mu,T)$ be a system, let
$A\subset X$ be a set with $\mu(A)>0$ and let $k\geq 1$ an integer.
$$
\liminf_{N\to+\infty}\frac 1{N}\sum_{n=0}^{N-1}
              \mu(A\cap T^nA\cap T^{2n}A\cap\dots\cap T^{(k-1)n}A)\ >0\ .
$$
\end{Fth}
In particular there exists $n\geq 1$ such that
$\mu(A\cap T^nA\cap T^{2n}A\cap\dots\cap T^{(k-1)n}A)>0$.
Then Furstenberg deduced Szemer\'edi's Theorem by using the following
Correspondence Principle (with  $m_j=(j-1)n$ for
$j=1,2,\dots,k$):

\begin{CorrPrinc}[\cite{F3}]
\label{prop:correspondence}
Let $E$ be a set of integers with positive upper density. There
exist a system $(X,\mu,T)$ and a subset $A$ of $X$ with
$\mu(A)= d^*(E)$ such that
$$
d^*\bigl((E+m_1)\cap\dots\cap(E+m_k)\bigr)
\geq
\mu( A\cap T^{m_1}A\cap\dots\cap T^{m_k}A)
$$
for all integers $k\geq 1$ and all $m_1,\dots,m_k\in\N$.
\end{CorrPrinc}

This indirect proof cannot provide  explicit bounds in the finite
version of Szemer\'edi's Theorem (Theorem~\ref{th:SZF}), but it is conceptually much simpler than the
original proof. Moreover the direction initiated by Furstenberg's 
Theorem led to several
generalizations, each of them inducing its combinatorial counterpart
by the Correspondence Principle (or by some variation of it).
For some some of these generalizations, there is still no known proof other that
the ergodic theoretic proof.
Below we only discuss 
two of these generalizations.

Recently Tao~\cite{Tao} gave a new proof of Szemer\'edi's Theorem that
can
be viewed as a cross between combinatorial and ergodic
methods: the proof is purely combinatorial in the sense that it deals
only with finite sets (subsets of $\Z/N\Z$). But the vocabulary and
the ``philosophy'' of the paper are much closer to ergodic theory.
The proof uses in particular an induction that mimics the inductive
construction of a sequence of extensions used in Furstenberg's paper. It
would be interesting to compare the way Tao uses Gowers' norms with
the way we use the ergodic seminorms in~\cite{HK4} (see
Section~\ref{sec:seminorms} and Appendix~\ref{subsec:gowersnorms}).

\subsubsection{The Polynomial Szemer\'edi Theorem}
In the next theorem proved by Bergelson and Leibman the exponents $n$, $2n$,
\dots, appearing in Furstenberg's Theorem are replaced by integer polynomials
$p_1(n)$, $p_2(n)$, \dots (an \emph{integer polynomial} is a
polynomial taking integer
values on the integers).

\begin{PMRT}[Bergelson \& Leibman~\cite{BL}]
Let $(X,\mu,T)$ be  system, let $\ell\geq 1$ be an integer and let 
$p_1(n)$, $p_2(n)$, \dots, $p_{\ell-1}(n)$ be integer polynomials
with $p_j(0) = 0$ for $j = 1, 2, \ldots, \ell-1$.
For every $A\subset X$ with $\mu(A)>0$ we have
$$
\liminf_{N\to\infty}\frac{1}{N}\sum_{n=0}^{N-1}\mu\bigl(A\cap T^{p_1(n)}A\cap
T^{p_2(n)}\cap\ldots\cap T^{p_{\ell-1}(n)}A \bigr) > 0 \ .
$$
\end{PMRT}

The Correspondence Principle immediately gives:

\begin{PST}[\cite{BL}]
Let $E$ be a set of integers with positive upper density and
$p_1(n)$, $p_2(n)$,\dots, $p_{\ell-1}(n)$ integer polynomials
with $p_j(0) = 0$ for $j = 1, 2, \ldots, \ell-1$. Then  there exist
integers $a$ and  $d>0$ such that
$$
    \{a,a+p_1(d),a+p_2(d),\dots a+p_{\ell-1}(d)\}\subset E\ .
$$
\end{PST}
Until now, the ergodic proof is the only known and in particular,
there is no known  combinatorial proof.

\subsubsection{The Multidimensional Szemer\'edi Theorem}
The following theorem of Furstenberg and Katznelson generalizes
Furstenberg's Theorem for
several commuting transformations.

\begin{MRT}[Furstenberg \& Katznelson~\cite{FK79}]
Let $k\geq 1$ be an integer and $T_1,T_2,\dots, T_k$ commuting
measure preserving
transformations of the probability space $(X,\mu)$. Then for any
subset $A$ of $X$ with $\mu(A)>0$ we have
$$
    \liminf_{N\to\infty}\frac 1N\sum_{n=0}^{N-1}
\mu( T_1^nA\cap T_2^nA\cap\dots\cap T_k^nA)>0\ .
$$
\end{MRT}
The original theorem corresponds to the case that
$T_j=T^{j-1}$ for
$j=1,\dots,k$.

The upper density of a subset of $\N^k$ is defined analogously
as for a subset of $\N$.
The combinatorial counterpart (see Note~\ref{note:Multi})
of the theorem above is:

\begin{MST}[\cite{FK79}]
Let $E\subset \N^k$ be a set of positive upper density. Then for any finite
subset $F$ of $\N^k$ there exists $\ba\in\N^k$ and an integer
$d>0$ such that $\ba+d.F\subset E$.
\end{MST}
Here, $\ba+d.F=\{\ba+d.\bx\colon\bx\in F\}$.
When $F=\{0,1,\dots,\ell-1\}^k$ the set $\ba+d.F$ can be called
an arithmetic progression of dimension $k$ and length $\ell$, hence
the name of the theorem. The first combinatorial proof of this result
was given by Gowers~\cite{Gowers2}.

There exist several other generalizations of Furstenberg's Theorem
and for each of them a generalization of Szemer\'edi's Theorem (see
for example~\cite{BMC}). The deepest result in this class is the
Density Hales-Jewett Theorem \cite{FK91} of Furstenberg and
Katznelson.

\subsection{Convergence results}
It is a natural question to ask whether the $\liminf$ in Furstenberg's
Theorem and in its generalizations are actually limits. Our main goal
in these lectures is to give an idea of the proof of the following
results, that we refer to as \emph{convergence theorems for multiple ergodic
averages}. The first theorem shows the convergence for arithmetic
progressions:
\begin{theorem}[Host \& Kra~\cite{HK4}]
\label{th:convprog}
Let $(X,\mu,T)$ be a system, $k\geq 1$ be an integer and
$f_1,f_2,\dots,f_k$ be bounded functions on $X$. Then the averages
\begin{equation}
\label{eq:AP}
\frac 1N\sum_{n=0}^{N-1}f_1(T^nx)f_2(T^{2n}x)\dots f_k(T^{kn}x)
\end{equation}
converge in $L^2(\mu)$.
\end{theorem}
Taking $f_1=f_2=\dots=f_k=\one_A$ and integrating~\eqref{eq:AP} over $A$
we get that the $\liminf$ in Furstenberg's Theorem is a limit.

The  convergence in $L^2(\mu)$ of these averages for $k=3$ with the added
hypothesis that the system is totally ergodic was shown by Conze and
Lesigne in  a series of
papers (\cite{CL1}, \cite{CL2}, \cite{CL3}, see also~\cite{Lesigne2})
and by Host and Kra~\cite{HK1} in the general case (see
also~\cite{FW}).
Furstenberg~\cite{F2} proved the convergence for every $k$ under the
assumption that the system is weakly mixing.

A similar convergence result also holds for polynomials:

\begin{theorem}[Host \& Kra~\cite{HK6}; Leibman~\cite{Lei3}]
\label{th:pol}
Let $(X,\mu,T)$ be a system, $k\geq 1$ an integer,
$p_1(n)$, $p_2(n)$,\dots, $p_k(n)$ integer polynomials
and $f_1,f_2,\dots,f_k$ bounded functions on $X$.
Then the averages
\begin{equation}
\label{eq:pol-averages}
\frac 1{N}\sum_{n=0}^{N-1}
f_1(T^{p_1(n)}x)f_2(T^{p_2(n)}x)\dots f_k(T^{p_k(n)}x)
\end{equation}
converge in $L^2(\mu)$.
\end{theorem}
It follows
that the averages appearing in the Polynomial Multiple
Recurrence Theorem converge.

The result of Theorem~\ref{th:pol} was proved by
Bergelson~\cite{Berg1} for weakly mixing
systems. Furstenberg and Weiss~\cite{FW} showed the convergence
for two particular cases when $k=2$: for $p_1(n)=n$, $p_2(n)=n^2$ and
for
$p_1(n)=n^2$, $p_2(n)=n^2+n$.  The paper~[HK6] contains a proof for the
general case, except when the system is not totally ergodic and at
least one
polynomial is of degree $1$ and some other is of degree $>1$. This
restriction was lifted by Leibman~\cite{Lei3}.
Frantzikinakis and Kra have shown that if the system is totally
ergodic and the polynomials $p_i$ are linearly independent, then the
limit of the
averages~\eqref{eq:pol-averages} is constant and equal to the product
of the integrals of the functions $f_i$.

The paper~\cite{HK4} also contains  the proof of convergence of another
type of multiple ergodic averages, the \emph{cubic averages}: see 
Appendix~\ref{note:cubic}.

In the above results the averages on $[0,N-1)$ can be replaced by
averages on any sequence of intervals whose lengths tend to infinity.
It is sufficient to prove these two theorems for ergodic systems, as
the general case follows by ergodic decomposition. So we henceforth
assume ergodicity.

The case of several commuting transformations remains essentially open.
The problem can be formulated as follows:
\begin{question}
Let $k\geq 1$ be an integer and $T_1,T_2,\dots, T_k$ commuting
measure preserving
transformations of the probability space $(X,\mu)$. Is it true that
for all bounded functions $f_1,f_2,\dots,f_k$ on $X$ the averages
$$
\frac 1N\sum_{n=0}^{N-1}
f_1(T_1^nx)f_2(T_2^nx)\dots f_k(T_k^nx)
$$
converge in $L^2(\mu)$?
\end{question}
The answer was shown to be positive for two transformations by Conze
and Lesigne~\cite{CL1}. Frantzikinakis and Kra~\cite{FK2} proved the
convergence for an arbitrary number of transformations under the
additional hypothesis that all the transformations $T_i$ and all the
transformations $T_iT_j\inv$ $i\neq j$, are ergodic. These  are very
strong hypotheses: it can be assumed without loss  that the
transformations are jointly ergodic but not that any individual
transformation is ergodic. The tools developed below for one
transformation do not generalize to the case of several 
transformations.

The convergence almost everywhere of the different averages considered
here is an open and probably very difficult problem. The
unique result in this direction is due to
Bourgain~\cite{Bourgain1}: the convergence a.e. of the
averages~\eqref{eq:AP} for $k=2$.

\subsection{Notes on Section~\ref{sec:context}}

\subsubsection{The analytic form}
\label{subsec:anal}
The finite version of Szemer\'edi's Theorem (Theorem~\ref{th:SZF}) follows easily from the
analytic form, but the converse implication is more tricky. The
following result is apparently stronger than Szemer\'edi's Theorem but
can be deduced from it (this is a nontrivial exercise):
\begin{theorem*}
For every integer $\ell\geq 2$ and every real $\delta>0$ there exists
a constant $c(\ell,\delta)>0$ such that:\\
For every integer $N\geq 1$, every subset $E$ of $\{1,2,\dots,N\}$
with more than $\delta N$ elements contains at least $c(\ell,\delta)N^2-1$
    arithmetic progression of length $\ell$.
\end{theorem*}
\noindent
(The $-1$ in the formula is simply a way to eliminate trivialities
for small $N$.)

This immediately implies the analytic form in the particular case that
the function $f$ takes its values in $\{0,1\}$. The general case
follows by standard methods.

\subsubsection{About the theorem of Green and Tao}
\label{subsec:noteGT}
The set  $P$ of primes satisfies
$$
   \sum_{p\in P}\frac 1p=+\infty
$$
and it is natural to ask wether every subset of $\N^*$ with a
divergent series of inverses contains arbitrarily long arithmetic
progressions. This question
was asked by Erd\"os and Tur\'an~\cite{ET} in 1936 and is still open.

\subsubsection{From the Multidimensional Recurrence Theorem to the
Multidimensional Szemer\'edi Theorem}
\label{note:Multi}
This implication uses a generalization of the Correspondence Principle
that we state here. Let $T_1,T_2,\dots, T_k$ be commuting
measure preserving
transformations of the probability space $(X,\mu)$. For
$\bn=(n_1,n_2,\dots,n_k)\in\Z^k$ we write
$T^\bn=T_1^{n_1}T_2^{n_2}\dots T_k^{n_k}$.
\begin{proposition*}
Let $E\subset\N^k$  be a set of positive upper density. There exist $k$ measure
preserving transformations $T_1,T_2,\dots,T_k$ of a probability
space $(X,\mu)$ and a subset $A$ of $X$ with $\mu(A)=d^*(E)$ and
$$
d^*\Bigl(\bigcap_{\bn\in F}(E+\bn)\Bigr)\geq
   \mu\Bigl(\bigcap_{\bn\in F}T^\bn A\Bigr)
\text{ for every finite subset $F$ of $\N^k$.}
$$
\end{proposition*}
Then apply the Recurrence Theorem with the commuting measure
preserving transformations $T^\bn$, $\bn\in F$.

\section{Nilmanifolds and nilsystems}
\label{sec:nil}
We present here a class of systems for which the convergence results
(Theorems~\ref{th:convprog} and ~\ref{th:pol}) can be proven in an easier
way. In the rest of these lectures we 
  explain without details how the general case
can be deduced from this particular one.

\subsection{Definitions and fundamental properties}
\label{subsec:defnil}
Let $G$ be a group. For $g,h\in G$ we write $[g,h]=g\inv h\inv gh$.
The \emph{lower central series}
$$
G=G_1\supset G_2\supset \dots\supset G_j\supset G_{j+1}\supset\dots
$$
of $G$ is defined by $G_1=G$ and, for $j\geq 1$,
$$
G_{j+1} \text{ is the subgroup of $G$ spanned by }
\{[g,h]\colon g\in G,\ h\in G_j\}\ .
$$
Let $k\geq 1$ be an integer. We say that \emph{$G$ is $k$-step
nilpotent} if $G_{k+1}=\{1\}$.

Let $G$ be a $k$-step nilpotent Lie group and let $\Lambda$ be a discrete
cocompact subgroup. The compact manifold $X=G/\Lambda$ is called a
\emph{$k$-step nilmanifold}.
The fundamental properties of nilmanifolds were established by
Malcev~\cite{Malcev}. We recall here only the properties we need.

  The group $G$ naturally acts on $X$ by
left translation and we write $(g,x)\mapsto g\cdot x$ for
this action.
There exists a unique  Borel probability measure $\mu$ on $X$ invariant
under this action, called the  \emph{Haar measure} of $X$.
We make use of the following property
which appears in~\cite{Malcev} for connected groups and is proved in
Leibman~\cite{Lei2} in a similar way for the general case:
\begin{itemize}
\item\em
For every integer $j\geq 1$, the subgroups $G_j$ and $\Lambda G_j$
are closed in
$G$. It follows that the group $\Lambda_j=\Lambda\cap G_j$ is
cocompact in $G_j$.
\end{itemize}

Let $t$ be a fixed element of $G$ and let $T:X\to X$ be the transformation
$x\mapsto t\cdot x$. Then $(X,T)$ is called a \emph{$k$-step topological
nilsystem} and $(X,\mu,T)$ is called a \emph{$k$-step nilsystem}.
All the notation introduced above is used freely throughout the rest
of this
section.

Fundamental properties of nilsystems were established by
Auslander, Green and
Hahn~\cite{AGH} and by Parry~\cite{Pa1}.
Further ergodic properties were proven by Parry~\cite{Pa2} and
Lesigne~\cite{Le} when the group $G$ is connected, and
generalized  by Leibman (\cite{Lei2}, \cite{Lei3}). We state the
results we need in the next two theorems.
\begin{theorem}
\label{th:nilsystems1}
The following properties are equivalent:
\begin{itemize}
\item
$(X,\mu,T)$ is ergodic.
\item
$(X,T)$ is uniquely ergodic, meaning that $\mu$ is the unique
$T$-invariant probability measure on $X$.
\item
$(X,T)$ is minimal, meaning that every orbit under $T$ is dense.
\item
$(X,T)$ is transitive, meaning that there exists at least one
dense orbit under $T$.
\end{itemize}
\end{theorem}
The second theorem  can be viewed as a
particular case of general results of Ratner (\cite{Ratner}) and
Shah (\cite{Sh}).
\begin{theorem}
\label{th:nilsystems2}
Let $x\in X$ and let $Y$ be the closure of the orbit of $x$ under $T$.
Then $(Y,T)$ can be given
the structure of a topological nilsystem, that is $Y=H/\Gamma$ where
$H$ is a closed Lie subgroup of $G$ containing $t$ and $\Gamma$ is a
closed cocompact subgroup of $H$.
\end{theorem}

By Theorem~\ref{th:nilsystems1}, $(Y,T)$ is uniquely ergodic. We
immediately deduce:
\begin{corollary}
\label{cor:nilsystems1}
For every continuous function $f$ on $X$ the averages
\begin{equation}
\label{eq:erg-nil}
\frac 1N\sum_{n=0}^{N-1} f(T^nx)
\end{equation}
converge for every $x\in X$ as $N\to +\infty$.
\end{corollary}

\subsection{Two examples}
We  review the simplest examples
of $2$-step nilsystems.
\subsubsection{}
\label{subsec:example1}
Let $G=\Z\times\T\times \T$,  with  multiplication given by
$$
(k,x,y)*(k',x',y')=(k+k', x+x', y+y'+ 2kx')\ .
$$
Then $G$ is a Lie group. Its commutator subgroup is
$\{0\}\times \{0\}\times\T$ and $G$ is $2$-step nilpotent.
The subgroup $\Lambda=\Z\times\{0\}\times\{0\}$ is discrete and
cocompact.
Let $X$ denote the nilmanifold $G/\Lambda$
and we maintain the notation of the preceding section.

Fix $\alpha\in\T$ and define  $t=(1,\alpha,\alpha)\in G$
and let $T:X\to X$ be the translation by $t$.
Then $(X,\mu,T)$ is a $2$-step nilsystem. It can be shown that
it is ergodic if and only if $\alpha$ is irrational.

We give an alternate description of this system.
The map $(k,x,y)\mapsto (x,y)$ from $G$ to $\T^2$ induces a
homeomorphism of $X$ onto $\T^2$.
Identifying $X$ with $\T^2$ via this homeomorphism, the measure
$\mu$ becomes equal to $m_\T\times m_\T$ where $m_\T$ is the Haar
measure of $\T$ and the transformation $T$ of $X$ is
given for $(x,y)\in\T^2=X$ by
$$
T(x,y)=(x+\alpha, y+2x+\alpha)\ .
$$

\subsubsection{}
\label{subsec:example2}
Let $G$ be the Heisenberg group $\R\times\R\times\R$, with
multiplication given by
\begin{equation}
\label{eq:heisen}
(x,y,z)*(x',y',z') = (x+x', y+y', z+z'+xy') \ .
\end{equation}
Then $G$ is a $2$-step nilpotent Lie group.
The subgroup $\Lambda = \Z\times\Z\times\Z$ is discrete and cocompact.
Let $X=G/\Lambda$ and let
  $T$ be the translation by some $t = (t_1, t_2, t_3)\in G$.
We have that $(G/\Lambda,T)$ is a nilsystem.
It can be showed that it is ergodic if and only if $t_1$ and $ t_2$
are independent over $\Q$.

\subsection{Convergence of multiple ergodic averages for nilsystems}

Let $(X,\mu,T)$ be a nilsystem. We use here the notation of
Subsection~\ref{subsec:defnil}.

Let $k\geq 2$ be an integer. Then $X^k$ can be given the structure of a
nilmanifold,  quotient of the Lie
group $G^k$ by its discrete cocompact subgroup $\Lambda^k$.
Set $s=(t,t^2,\dots,t^k)\in G^k$ and let $S$ be the translation by $s$
on $X^k$. Then $(X^k,S)$ is a topological nilsystem.

Let $f_1,f_2,\dots,f_k$ be continuous functions on $X$.
By applying   Theorem~\ref{th:nilsystems2} to the nilsystem
$(X^k ,S)$ and to the continuous function
$(x_1,x_2,\dots,x_k)\mapsto f_1(x_1)f_2(x_2)\dots f_k(x_k)$ at the
point $(x,x,\dots,x)\in X^k$ we get that the average
$$
\frac 1N\sum_{n=0}^{N-1}f_1(T^nx)f_2(T^{2n}x)\dots f_k(T^{kn}x)
$$
converges for every $x\in X$.
By a standard density argument we have:
\begin{corollary}
\label{cor:convprognil}
Theorem~\ref{th:convprog} holds for  nilsystems.
\end{corollary}

An explicit expression of the limit was given by
Ziegler~\cite{Ziegler} (a shorter proof can be found in~\cite{BHK}).

  A similar result holds for the polynomial averages:
\begin{corollary}[\cite{Lei2}]
\label{cor:convpolnil}
Theorem~\ref{th:pol} holds for nilsystems.
\end{corollary}
In place of  Corollary~\ref{cor:nilsystems1}, the  proof uses an
extension due to Leibman of this result for \emph{polynomial sequences} in a
nilmanifold.

\section{Some seminorms}
\label{sec:seminorms}
In this Section and the next ones $(X,\mu,T)$ is an ergodic
system.
We introduce a sequence of seminorms of
$L^\infty(\mu)$ that we use to bound the different averages under
consideration. These seminorms should be compared with the norms
introduced by Gowers in a completely different context:
see Appendix~\ref{subsec:gowersnorms}.

\subsection{Notation}
We need some
notation to be used throughout the remainder of these notes.

We write $C\colon\C\to\C$ for the conjugacy map $z\mapsto\bar z$.
Let $k\geq 1$ be an integer. We write
$\eps=\eps_1\eps_2\dots\eps_k$ with
$\eps_i\in\{0,1\}$ for a point of $\{0,1\}^k$, without commas or
parentheses and $\lvert\eps\rvert=\eps_1+\eps_2+\dots+\eps_k$.
$\zero$ denotes the point $00\dots 0\in\{0,1\}^k$.

For every integer $k\geq 0$ we write $X\type k=X^{2^k}$ and for $k\geq
1$ the points
of $X\type k$ are written
$\bx=(x_\eps\colon\eps\in\{0,1\}^k)$. We write $T\type k$
for the transformation $T\times T\times\dots\times T$ ($2^k$ times) of
this space. We often identify $X\type{k+1}$ with $X\type k\times X\type
k$ writing $\bx=(\bx',\bx'')$ for a point of $X\type {k+1}$ where
$\bx',\bx''\in X\type k$ are defined by
$$
\text{for every }\eps\in\{0,1\}^k,\
x'_\eps=x_{\eps 0}\text{ and } x''_\eps =
x_{\eps1}\ .
$$

\subsection{Construction of some measures and some seminorms}
\label{subsec:constsemi}

We define by induction a $T\type k$-invariant measure $\mu\type k$ on
$X\type k$ for every integer $k\geq 0$.
Set $\mu\type 0=\mu$. Assume that $\mu\type k$ is defined for some
$k\geq 0$.
Let $\CI\type k$ denote the
$T\type k$ invariant $\sigma$-algebra of
$(X\type k, \mu\type k,T\type k)$.
Identifying $X\type{k+1}$ with $X\type k\times X\type k$ as
explained above, we define the system
$(X\type{k+1},\mu\type{k+1},T\type{k+1})$ to be the \emph{relatively
independent joining} of two copies of $(X\type k, \mu\type k,T\type k)$
over $\CI\type k$.

This means that $\mu\type{k+1}$ is the measure on
$X\type{k+1}=X\type k\times X\type k$ characterized by:

\emph{when $F'$ and $F''$ are bounded
functions on $X\type k$,}
\begin{equation}
\label{eq:defmes}
\int_{X\type{k+1}} F'(\bx')F''(\bx'')\,d\mu\type{k+1}(\bx)=
\int_{X\type k}\E(F'\mid\CI\type k)\,\E(F''\mid\CI\type k)\,d\mu\type
k\ .
\end{equation}
This measure  is invariant under
$T\type{k+1}=T\type k\times T\type k$
and each of its two natural projections on $X\type k$
is equal to $\mu\type k$.
Note that when $F$ is a function on $X\type k$, measurable with
respect to $\CI\type k$ that is invariant under
$T\type k$, we have
\begin{equation}
\label{eq:equalinvar}
F(\bx')=F(\bx'')\text{ for $\mu\type{k+1}$-almost every
}\bx=(\bx',\bx'')\in X\type{k+1}\ .
\end{equation}

By induction, each of the $2^k$ natural
  projections of $\mu\type k$ on $X$  is equal to $\mu$.
It follows immediately from this definition that  for every bounded
function $f$ on $X$ the integral
$$
   \int_{X\type k}\prod_{\eps\in\{0,1\}^k}
C^{\lvert\eps\rvert}f(x_\eps)\,d\mu\type k(\bx)
$$
is real and nonnegative and we can define
\begin{equation}
\label{eq:defsemi}
\nnorm f_k = \Bigl(\int_{X\type k}\prod_{\eps\in\{0,1\}^k}
C^{\lvert\eps\rvert}f(x_\eps)\,d\mu\type k(\bx)\Bigr)^{1/2^k}\ .
\end{equation}

As $X$ is assumed to be ergodic, the $\sigma$-algebra $\CI\type 0$ is
trivial and $\mu\type 1=\mu\times\mu$. We therefore have
$$
   \nnorm f_1=\Bigl(\int_{X\times X}f(x_0)\overline{f(x_1)}
   \,d\mu\times\mu(x_0,x_1)
\Bigr)^{1/2}=\Bigl|\int f(x)\,d\mu(x)\Bigr|\ .
$$

  From the inductive definition~\eqref{eq:defmes} of the measures
we get that
$$
\nnorm f_{k+1}\geq\nnorm f _k \text{ for every $f$ and every $k$.}
$$
The following Lemma follows immediately from the definition of the
measures and the Ergodic Theorem.
\begin{lemma}
\label{lemma:recur-mes}
For every integer $k\geq 0$ and every bounded function $f$ on $X$,
\begin{equation}
\label{eq:recur-mes}
\nnorm f_{k+1}=\Bigl(\lim_{N\to\infty} \sum_{n=0}^{N-1}
\nnorm{f\,.\,T^n\overline{f}\,}^{2^k}_k\Bigr)^{1/2^{k+1}}\ .
\end{equation}
\end{lemma}

\subsection{The Kronecker factor and the measure
   $\mu\type 2$}

\subsubsection{The notion of a factor}
As usual in ergodic theory, we use the term \emph{factor} with two
different but equivalent meanings. First, a factor of the system
$(X,\mu,T)$ is a $T$-invariant sub-$\sigma$-algebra $\CY$ of $\CX$.

On the other hand, if $(Y,\nu,S)$ is a system and $\pi\colon X\to Y$
is
a measurable map mapping $\mu$ to $\nu$ and such that $S\circ\pi=\pi\circ T$
($\mu$-a.e.),
then $\pi$ is called a \emph{factor map} and $Y$ is also called a
factor of $X$.
In this situation we always identify the $\sigma$-algebra $\CY$ of $Y$
with its inverse image $\pi\inv(\CY)$, which is an invariant
sub-$\sigma$-algebra of $\CX$, that is a factor of $X$ under the first
definition.
It is thus natural to denote the transformation on $Y$ by the same
letter as the
transformation on $X$, meaning by $T$ in our case.

It can be shown that every invariant sub-$\sigma$-algebra of $\CX$
can be associated to a factor map in this way and thus the two
definitions are functionally equivalent. We pass freely from one
definition to the other.

Let $f$ be an integrable function on $X$. We consider $\E(f\mid\CY)$
as a function defined on $X$ and
we write $\E(f\mid Y)$ for the function on $Y$ defined by
$\E(f\mid Y)\circ\pi=\E(f\mid \CY)$.
It is characterized by:
$$
  \forall g\in L^\infty(\nu),\
\int_Y \E(f\mid Y)(y). g(y)\,d\nu(y)=\int_X f(x).g(\pi(x))\,d\mu(x)\ .
$$
\subsubsection{The Kronecker factor}
\label{subsec:kro}
The Kronecker factor of the system $X$ is written $Z_1(X)$ or
$Z_1$ when the system under consideration is clear from the context.
We recall here the definition and some classical properties.

Viewed as a $\sigma$-algebra, the \emph{Kronecker factor} $\CZ_1$ of
$\CX$ is defined to be the sub-$\sigma$-algebra of $\CX$ generated by the
eigenfunctions of this system; is is also the smallest
sub-$\sigma$-algebra of $\CX$ such that all the invariant functions
of the system $(X\times X,\mu\times\mu,T\times T) $ are measurable with
respect to $\CZ_1\otimes\CZ_1$.

When considered as a system, the Kronecker factor $(Z_1,m,T)$
is a \emph{rotation}; this means that $Z_1$ is a compact abelian
group  with Haar measure $m$ and
that the transformation $T$ has the form $z\mapsto z+\alpha$  where
$\alpha$ is a fixed element of $Z_1$.

We write $\pi_1\colon X\to Z_1$
for the factor map. Then every eigenfunction of $X$ has the form
$f(x)=c.\chi(\pi_1(x))$ where $c$ is a constant and $\chi$ is a
\emph{character} of $Z_1$, that is a continuous group homomorphism
from $Z_1$ to the circle group $\CS^1$. Every $T\times T$-invariant
function on $X\times X$ can be written
$f(x_0,x_1)=g(\pi(x_1)-\pi(x_0))$ where $g$ is a function on $Z_1$.
Therefore, when $f_0$ and $f_1$ are bounded functions on $X$,
\begin{multline*}
   \E(f_0\otimes f_1\mid\CI\type 1)(x_0,x_1)\\
=\int_{Z_1}\E(f_0\mid Z_1)(z)\,
\E\bigl(f_1\mid Z_1)(z+\pi(x_1)-\pi(x_0)\bigr)\,dm(z)\ .
\end{multline*}
\subsubsection{The measure $\mu\type 2$}
\label{subsec:mu2}
We deduce a more explicit expression for the measure $\mu\type 2$:
when $f_\eps$, $\eps\in\{0,1\}^2$ are four measurable
functions on $X$,
\begin{multline}
\label{eq:mu2}
\int_{X\type 2}f_{00}(x_{00})f_{01}(x_{01})f_{10}(x_{10})f_{11}(x_{11})
\,d\mu\type 2(\bx)\\
=\int_{Z_1\times Z_1\times Z_1}
\tilde f_{00}(z)\,
\tilde f_{01}(z+s)\,
\tilde f_{10}(z+t)\,
\tilde f_{11}(z+s+t)\,
dm(z)\,dm(s)\,dm(t)
\end{multline}
where we write $\tilde f_\eps$ for $\E(f_\eps\mid Z_1)$.
This implies that $\nnorm f_2$ is the $\ell^4$-norm of the Fourier Transform of
$\E(f\mid Z_1)$. In particular:
\begin{lemma}
\label{lemma:mu2}
For every bounded function $f$ on $X$,
$\nnorm f_2=0$ if and only if
$\E(f\mid Z_1)=0$.
\end{lemma}

We give one more formula for the measure $\mu\type 2$.
   For every $s\in Z_1$ we define a probability measure $\mu_s$ on
$X\times X$ by
\begin{equation}
\label{eq:def-mus}
\int_{X\times X} f_0(x_0)f_1(x_1)\,d\mu_s(x_0,x_1)=\int_{Z_1}
\E(f_0\mid Z_1)(z)\,\E(f_1\mid Z_1)(z+s)\,dm(z)\ .
\end{equation}
This measure is invariant under $T\times T$ and we have
\begin{equation}
\label{eq:ergod-dec}
\mu\times\mu=\int_{Z_1} \mu_s\,dm(s)\ .
\end{equation}
  From the remarks above it follows that this formula is the
\emph{ergodic decomposition} of $\mu\times\mu$ under $T\times T$; in
particular the system $(X\times X,\mu_s,T\times T)$ is ergodic for
$m$-almost every $s\in Z_1$. We have:
$$
   \mu\type 2=\int_{Z_1}\mu_s\times\mu_s\,dm(s)
$$
(notice that $\mu_s\times\mu_s$ is a measure on $(X\times X)\times
(X\times X)=X\type 2$).

\subsection{Seminorms and arithmetic progressions}

Later we prove  that for every $k\geq 0$ the map $f\mapsto\nnorm
f_k$ is a seminorm on $L^\infty(\mu)$. Admitting this for the moment
we show now how these seminorms arise in the questions of
convergence we are studying.

\begin{proposition}
\label{prop:semiprog}
Let $f_1,f_2,\dots,f_k$ be bounded functions on $X$ with
$\norm{f_\ell}_\infty\leq 1$ for $\ell=1,2,\dots,k$. Then
\begin{equation}
\label{eq:semiprog}
\limsup_{N\to+\infty}
\Bigl\Vert
\frac 1N\sum_{n=0}^{N-1}f_1(T^nx)f_2(T^{2n}x)\dots f_k(T^{kn}x)
\Bigr\Vert_{2}
\leq \min_{1\leq \ell\leq k}\bigl(\ell.\nnorm{f_\ell}_k\bigr)\ .
\end{equation}
\end{proposition}

The proof relies on an iterated use of a Hilbert space variant of the
van der Corput Lemma.

\begin{VdC}[\cite{Berg1}]
Let $\{\xi_n\}$ be  a sequence in a Hilbert space $\CH$, with
$\norm{\xi_n}\leq 1$ for every $n$. Then
\begin{equation}
\label{eq:VdC}
   \limsup_{N\to+\infty}\Bigl\Vert\frac
   1N\sum_{n=0}^{N-1}\xi_n\Bigr\Vert^2
\leq \limsup_{H\to+\infty}\frac 1H\sum_{h=0}^{H-1}
\limsup_{N\to+\infty}\Bigl|\frac 1N\sum_{n=0}^{N-1} \langle
\xi_{n}\mid\xi_{n+h}\rangle \Bigr|
\end{equation}
\end{VdC}
\begin{proof}[Proof of Proposition~\ref{prop:semiprog}]
We proceed by induction. For $k=1$ the bound is given by the Ergodic
Theorem and the definition of $\nnorm\cdot_1$. Let $k\geq 2$ and
assume that the bound holds for $k-1$.
Let $f_1,\dots,f_k$ be as in the theorem, and choose
$\ell\in\{2,\dots,k\}$ (the case $\ell=1$ is similar). Write
$$
   \xi_n(x)=f_1(T^nx)f_2(T^{2n}x)\dots f_k(T^{kn}x)\ .
$$
For every $h\geq 0$, by using the Cauchy-Schwarz Inequality and the invariance
of $\mu$ under $T^n$ we get
$$
\Bigl|\frac 1N\sum_{n=0}^{N-1} \langle
\xi_{n}\mid\xi_{n+h}\rangle \Bigr|
\leq
\norm{ f_1.\overline{f_1}\circ T^h}_2\,.\,
\Bigl\Vert \sum_{n=0}^{N-1}
\Bigl(\prod_{i=2}^k(f_i.\overline{f_i}\circ T^{ih})\circ T^{(i-1)n}\Bigr)
\Bigr\Vert_2\ .
$$
By the inductive assumption,
$$
   \limsup_{N\to+\infty}\Bigl|\frac 1N\sum_{n=0}^{N-1} \langle
\xi_{n}\mid\xi_{n+h}\rangle \Bigr|
\leq \ell.
\nnorm{f_\ell. \overline{f_\ell}\circ T^{\ell h}}_{k-1}\ .
$$
By the van der Corput  Lemma,
\begin{multline*}
\limsup_{N\to+\infty}\Bigl\Vert\frac
   1N\sum_{n=0}^{N-1}\xi_n\Bigr\Vert_2^2
\leq
\ell. \limsup_{H\to+\infty}\frac 1H\sum_{h=0}^{H-1}
\nnorm{f_\ell. \overline{f_\ell}\circ T^{\ell h}}_{k-1}\\
\leq \ell^2. \limsup_{H\to+\infty}\frac 1H\sum_{h=0}^{H-1}
\nnorm{f_\ell. \overline{f_\ell}\circ T^{ h}}_{k-1}\\
\leq\ell^2.\limsup_{H\to+\infty}\Bigl(\frac 1H\sum_{h=0}^{H-1}
\nnorm{f_\ell. \overline{f_\ell}\circ T^{ h}}_{k-1}^{2^{k-1}}\Bigr)^{1/2^{k-1}}
\end{multline*}
and this last term is equal to $\ell^2 .\nnorm{f_\ell}_k^2$ by
Lemma~\ref{lemma:recur-mes}.
\end{proof}

\subsection{Seminorms and polynomial averages}
A similar bound holds for the averages~\eqref{eq:pol-averages}
considered in Theorem~\ref{th:pol}.
\begin{proposition}
\label{prop:semipol} Let $k\geq 1$ be an integer and
$p_1,p_2,\dots,p_k$ be integer nonconstant
polynomials such that for every $1\leq i\neq j\leq k$ the polynomial
$p_i-p_j$ is not constant. There exists  an integer $\ell\geq 1$ such
that for any bounded functions $f_1,f_2,\dots,f_k$ on $X$,
$$
\frac 1{N}\sum_{n=0}^{N-1}
f_1(T^{p_1(n)}x)f_2(T^{p_2(n)}x)\dots f_k(T^{p_k(n)}x)
$$
converges to zero in $L^2(\mu)$ whenever $\nnorm {f_i}_\ell=0$ for at
least one value of $i\in\{1,2,\dots,k\}$.
\end{proposition}
The proof uses induction on the \emph{polynomial family}
$(p_1,p_2,\dots,p_k)$, by using the van der Corput Lemma and the
Cauchy-Schwarz Inequality at each step. However this induction is much
more intricate than for arithmetic progressions and we do not present
it here.

\subsection{Invariance properties}
Let $\{0,1\}^k$  be identified with the set of  vertices of
the unit cube $[0,1]^k$.

   We call \emph{the group of the $k$-cube} the group of isometries of the
Euclidean cube $[0,1]^k$. We consider this group as acting on $X\type
k$ in the following way. Each element $\sigma$ of the group induces
a permutation, written $\sigma$ also, of the set $\{0,1\}^k$ of
vertices and this permutation in turn induces a transformation
$\sigma_*$ of $X\type k$ by:
$$
\text{for every }\eps\in\{0,1\}^k,\
   \bigl(\sigma_*\bx)_\eps=x_{\sigma(\eps)}\ .
$$

\begin{lemma}
\label{lemma:symm}
For every integer $k$ the measure $\mu\type k$ is invariant under the action
of the group of the $k$-cube.
\end{lemma}
An inequality similar to the classical Cauchy-Schwarz inequality
can be proven inductively by using this symmetry property:
\begin{lemma}
\label{lemma:CS}
   If $f_\eps$, $\eps\in\{0,1\}^k$, are $2^k$ bounded
functions on $X$ then
\begin{equation}
\label{eq:CS}
\Bigl|\int\prod_{\eps\in\{0,1\}^k}f_\eps(x_\eps)\,d\mu\type
k(\bx)\Bigr|\leq
\prod_{\eps\in\{0,1\}^k}\nnorm{f_\eps}_k\ .
\end{equation}
\end{lemma}
   It follows that the map $f\mapsto\nnorm f_k$ is subadditive:
\begin{corollary*}
For every $k\geq 1$, $\nnorm\cdot_k$ is a seminorm on $L^\infty(\mu)$.
\end{corollary*}

We use the geometric vocabulary for subsets of $\{0,1\}^k$;
for example, a \emph{side} is a subset of $\{0,1\}^k$ of the form
$\{\eps\colon\eps_i=0\}$ or $\{\eps\colon\eps_i=1\}$
for some $i\in\{1,\dots,k\}$. When $\alpha$ is a side of $\{0,1\}^k$
we define the \emph{side transformation} $T\type k _\alpha$ of $X\type
k$ by:
\begin{equation}
\label{eq:side}
\bigl(T\type k_\alpha\bx\bigr)_\eps=
\begin{cases}
Tx_\eps & \text{ if } \eps\in\alpha;\\
x_\eps & \text{ otherwise.}
\end{cases}
\end{equation}
Notice that the product of two side transformations corresponding to
opposite sides is equal to $T\type k$.

When $\alpha$ is the side $\{0,1\}^{k-1}\times\{1\}$,
$T\type k_\alpha$ is equal to $\id\type{k-1}\times T\type{k-1}$ under
the identification of $X\type k$ with $X\type{k-1}\times X\type{k-1}$.
By definition of the measure $\mu\type k$, it is invariant under this
transformation. As the group of the $k$-cubes leaves the measure
$\mu\type k$ invariant and acts transitively on the sides, we get:
\begin{lemma}
\label{lemma:side}
For every integer $k$ the  measure $\mu\type k$ is invariant under all side
transformations.
\end{lemma}
By induction, it can be checked that:
\begin{lemma}
\label{lemma:ergod}
The measure $\mu\type k$ is ergodic under the joint action of the side
transformations.
\end{lemma}

\subsection{Notes on Section~\ref{sec:seminorms}}
\label{subsec:notemeasurenil}

The measures $\mu\type k$ can be described explicitly in the case
that the given ergodic system $(X,\mu,T)$ is a $k$-step nilsystem.
We put $X=G/\Lambda$ and use the notation of Section~\ref{sec:nil}.

Let $k\geq 1$ be an integer. For every $g\in G$ and every side
$\alpha$ of $\{0,1\}^k$ we define an element $g\type k_\alpha$ of
$G\type k$ by
$$
(g\type k_\alpha)_\eps=\begin{cases}
g & \text{ if }\eps\in\alpha\ ;\\
1 & \text{ otherwise.}
\end{cases}
$$
We define the \emph{side group} of dimension $k$ to be the subgroup
$G\type k_{k-1}$ of $G\type k$ spanned by all
these elements. It can be checked that $G\type k_{k-1}$ is a closed
Lie subgroup of $G\type k$ and that the subgroup $\Lambda\type
k_{k-1}$ of $\Lambda\type k$ defined in the same way is equal to
$\Lambda\type k\cap G\type k_{k-1}$ and is discrete and cocompact in
$G\type k_{k-1}$. The nilmanifold $X_k=G\type k_{k-1}/\Lambda\type
k_{k-1}$ is naturally imbedded in $X\type k$ and the measure
$\mu\type k$ is the Haar measure of this nilmanifold.

\section{Building factors}
\label{sec:factors}
In this Section we continue to assume that $(X,\mu,T)$ is an ergodic
system. Lemma~\ref{lemma:mu2} gives a simple relation between the
seminorm $\nnorm\cdot_2$ and the Kronecker factor  $Z_1$ introduced
in Subsection~\ref{subsec:kro}.
For every $k\geq 2$ we define here a factor $Z_k$ of $X$ with a
similar relation to the seminorm $\nnorm \cdot_{k+1}$.

\subsection{The definition of the factors $Z_k$}
Let $k\geq 1$ be an integer. In the following construction the
coordinate of $X\type k$ indexed by $\zero$ plays a particular role.
Due to the symmetry of the measure $\mu\type k$, any other choice 
is possible with the same results, up to obvious changes in notation.

We note $X\typ k=X^{2^k-1}$ and write  each point $\bx\in X\type k$
as $\bx=(x_\zero,\tilde \bx)$, where
$$
\tilde\bx=
(x_\eps\colon\eps\in\{0,1\}^k\setminus\{\zero\})\in X\typ k
$$
and thus we identify $X\type k$ with $X\times X\typ k$.

Recall that
the projection of $\mu\type k$ on $X$ is equal to $\mu$.
We write $\mu\typ k$ for the projection of  $\mu\type k$
on $X\typ k$. This measure is invariant under $T\typ k=T\times
T\times\dots\times T$ ($2^k-1$ times). We say that the system
$(X\type k,\mu\type k, T\type k)$ is a \emph{joining} of the systems
$(X,\mu,T)$ and $(X\typ k,\mu\typ k,T\typ k)$.

Let $T\type k_1,T\type k_2,\dots,T\type k_k$ be the side transformations
of $X\type k$ corresponding to the $k$ sides of $\{0, 1\}^k$ not containing
$\zero=(0,0,\dots,0)$:
For $1\leq i \leq k$, $\bx\in X\type k$ and $\eps\in\{0,1\}^k$,
$$
  (T\type k_i\bx)_\eps=\begin{cases}
x_\eps & \text{ if }\eps_i=0\ ;\\
Tx_\eps & \text{ if }\eps_i=1\ .
\end{cases}
$$

For $1\leq i\leq k $ the transformation $T\type k_i$ leaves the
coordinate indexed by $\zero$ invariant and thus it can be written
$\id\times T\typ k_i$ for some measure preserving transformation
$T\typ k_i$ of $X\typ k$. Let $\CJ\typ k$ be the $\sigma$-algebra of
subsets of $X\typ k$ which are invariant under the transformations
$T\typ k_i$, $1\leq i\leq k$.
An induction using  relation~\eqref{eq:equalinvar}
gives:
\begin{lemma}
\label{lemma:Jk}
Let $B$ be a subset of $X\typ k$. Then $B$ belongs to $\CJ\typ k$ if
and only if there exists a subset $A$ of $X$ with
\begin{equation}
\label{eq:defZk}
\one_A(x_\zero)=\one_B(\tilde\bx)\text{ for $\mu\type k$-almost every
$\bx\in X\type k$.}
\end{equation}
\end{lemma}
This relation between $B$ and $A$ defines a bijection
(up to null sets) between the $\sigma$-algebra $\CJ\typ k$ and some
$\sigma$-algebra of $X$. We define:
\begin{definition}
$\CZ_{k-1}(X)$ is the sigma-algebra of subsets $A$ of $X$ such that
equality~\eqref{eq:defZk} holds for some subset $B$ of $X\typ k$.
\end{definition}
$\CZ_1(X)$ was already defined to be the Kronecker factor of $X$.
Below (Corollary~\ref{cor:Zk}) we show that the two definitions
coincide.
We write $\CZ_{k-1}$ instead of $\CZ_{k-1}(X)$ whenever it is possible
without ambiguity.
This $\sigma$-algebra is clearly invariant under $T$ and so it is
a factor map.
Let the  associated factor system be denoted by $(Z_{k-1}(X),\mu_{k-1},T)$ or
by $(Z_{k-1},\mu_{k-1},T)$
and let $\pi_{k-1}\colon X\to Z_{k-1}$ be the factor map.

\subsection{Elementary properties}
\begin{proposition}
\label{prop:propertiesZk}
\strut
\begin{enumerate}
\item\label{it:Zk1}
Consider   the $\sigma$-algebra $\CZ_{k-1}$ on $X$  and the
$\sigma$-algebra $\CJ\typ k$ on $X\typ k$ are identified by 
relation~\eqref{eq:defZk}. Then $(X\type k,\mu\type k)$ is the
relatively independent joining of $(X,\mu)$ and $(X\typ k,\mu\typ k)$
over their common $\sigma$-algebra $\CZ_{k-1}=\CJ\typ k$.
\item\label{it:Zk2}
For a bounded function $f$ on $X$, $\nnorm f_k=0$ if and only if
$\E(f\mid\CZ_{k-1})=0$.
\item\label{it:Zk3}
The measure $\mu\type k$ is relatively independent over its
projection on $Z_{k-1}\type k$; this means that when $f_\eps$,
$\eps\in\{0,1\}^k$, are bounded functions on $X$ then
\begin{equation}
\label{eq:Zk3}
\int_{X\type k}\prod_{\eps\in\{0,1\}^k}f_\eps(x_\eps)\,
d\mu\type k(\bx)=
\int_{X\type
k}\prod_{\eps\in\{0,1\}^k}\E(f_\eps\mid\CZ_{k-1})(x_\eps)\,
d\mu\type k(\bx)\ .
\end{equation}
Moreover, $\CZ_{k-1}$ is the smallest factor of $X$ with
this property.
\item\label{it:Zk4}
Every $T\type {k-1}$-invariant subset of $X\type {k-1}$ is measurable 
with  respect
to $\CZ_{k-1}\type {k-1}$. Moreover $\CZ\type {k-1}$ is the smallest
factor of $X$ with this property.
\end{enumerate}
\end{proposition}

\begin{proof}\strut

\ref{it:Zk1}
The meaning of this statement is perhaps not obvious and we begin with
some explanation. We have already introduced the factor map
$\pi_{k-1}\colon X\to Z_{k-1}$. As we identify $\CZ_{k-1}$ with the
$\sigma$-algebra $\CJ\typ k$ on $X\typ k$ we have also a factor map
$p_{k-1}\colon X\typ k\to Z_{k-1}$ with
$\CJ\typ k=p_{k-1}\inv(\CZ_{k-1})$. Relation~\eqref{eq:defZk}
can be written
\begin{equation}
\label{eq:defZkbis}
\pi_{k-1}(x_\zero)=p_{k-1}(\tilde\bx)\text{ for $\mu\type k$-almost
every }\bx=(x_\zero,\tilde\bx)\ .
\end{equation}
When $F$ is an integrable function on $X\typ k$, following our standard
convention we write $\E(F\mid Z_{k-1})$ for the function on $Z_{k-1}$
given by $\E(F\mid Z_{k-1})\circ p_{k-1}=\E(F\mid \CJ\typ k)$.
Statement~\ref{it:Zk1} means that for every  bounded function $f$
on $X$ and every bounded function $F$ on $X\typ k$,
\begin{equation}
\label{eq:intfF}
\int f(x_\zero)\,F(\tilde\bx)\,d\mu\type k(\bx)=
\int \E(f\mid Z_{k-1})\,\E(F\mid Z_{k-1})\,d\mu_{k-1}\ .
\end{equation}
This relation is similar to formula~\eqref{eq:defmes} used to
define the measure $\mu\type{k+1}$ in Section~\ref{subsec:constsemi}.

Let $f$ and $F$ be as above.
For $i=1,2,\dots,k$ the measure $\mu\type k$ and the function
$\bx\mapsto f(x_\zero)$ on $X\type k$ are  invariant under the
transformation $T\type k_i=\id\times T\typ k_i$. Therefore the first
integral in~\eqref{eq:intfF} remains unchanged when the function $F$
is replaced by its conditional expectation with respect to $\CJ\typ
k$ and this integral is equal to
$$
\int f(x_\zero)\,\E(F\mid Z_{k-1})\circ p_{k-1}(\tilde\bx)\,d\mu\type
k(\bx)\ .
$$
By using~\eqref{eq:defZkbis} we rewrite this integral as
$$
\int f(x_\zero)\,\E(F\mid Z_{k-1})\circ \pi_{k-1}(x_\zero)\,
d\mu(x_\zero)\ .
$$
In this last integral we can replace the function $f$ by its
conditional expectation with respect to $\CZ_{k-1}$ and 
equality~\eqref{eq:intfF} follows.

\ref{it:Zk2}
Assume that $\E(f\mid\CZ_{k-1})=0$. Using~\eqref{eq:intfF} with
$F(\tilde\bx)$ equal to the product of the functions
$C^{\lvert\eps\rvert}f(x_\eps)$
for $\eps\in\{0,1\}^k\setminus\{\zero\}$, we get that $\nnorm f_k=0$.

Assume conversely that $\nnorm f_k=0$. By using Lemma~\ref{lemma:CS}
and a density argument, we get that the
first integral in~\eqref{eq:intfF} is equal to zero for any choice
of the function $F$ on $X\typ k$ and thus in particular when
$F=\E(f\mid Z_{k-1})\circ p_{k-1}$.  In this case we have
$\E(F\mid Z_{k-1})=\E(f\mid Z_{k-1})$ and the second integral
in~\eqref{eq:intfF} is equal to  $\int \E(f\mid Z_{k-1})^2\,d\mu_{k-1}$.
As it is also equal to zero we have $\E(f\mid Z_{k-1})=0$ and thus
$\E(f\mid\CZ_{k-1})=0$.

\ref{it:Zk3}
By again using  the equality~\eqref{eq:intfF} we get that the first
integral in~\eqref{eq:Zk3} remains unchanged when the function
$\E(f_\zero\mid\CZ_{k-1})$ is substituted for the function
$f_\zero$. The same properties holds for the other vertices
$\eps\in\{0,1\}^k$ because of the symmetry of the measure
$\mu\type k$ (see Lemma~\ref{lemma:symm}) and this gives~\eqref{eq:Zk3}.

Let $\CY$ be a factor of $X$ with the same property. For each bounded
function on $X$ with $\E(f\mid \CY)=0$, equality~\eqref{eq:Zk3}
with all functions equal to $f$ gives  $\nnorm f_k=0$, and thus $\E(f\mid
\CZ_{k-1})=0$ by~\ref{it:Zk2}. This shows that $\CY\supset\CZ_{k-1}$
and achieves the proof of~\ref{it:Zk3}.

\ref{it:Zk4} Let $A$ be a $T\type{k-1}$-invariant subset of
$X\type{k-1}$. By construction of the measure $\mu\type k$ we have
$\one_A(\bx')=\one_A(\bx'')$ for $\mu\type k$ almost every
$\bx=(\bx',\bx'')\in X\type k$ and thus $\mu\type k(A\times
A)=\mu\type{k-1}(A)$. On the other hand,
\begin{multline*}
\mu\type k(A\times A)=\int \one_A\otimes\one_A\,d\mu\type k=
\int \E(\one_A\otimes\one_A\mid\CZ_{k-1}\type k)\,d\mu\type k\\
=\int \E(\one_A\mid\CZ_{k-1}\type{k-1})\otimes
\E(\one_A\mid\CZ_{k-1}\type{k-1})\, d\mu\type k
\end{multline*}
where the second equality holds because the measure $\mu\type k$ is
relatively independent over $Z_{k-1}\type k$. Moreover, as the
$\sigma$-algebra $\CZ_{k-1}\type{k-1}$ and the set $A$ are invariant
under $T\type{k-1}$, the function $\E(\one_A\mid\CZ_{k-1}\type{k-1})$
is invariant under this transformation and by construction of the measure
$\mu\type k$ we have $\E(\one_A\mid\CZ_{k-1}\type{k-1})(\bx')=
\E(\one_A\mid\CZ_{k-1}\type{k-1})(\bx'')$ for $\mu\type k$ almost every
$\bx=(\bx',\bx'')\in X\type k$. The last integral is therefore equal
to $\int \E(\one_A\mid\CZ_{k-1}\type{k-1})^2\,d\mu\type{k-1}$.
We get that the function $\one_A$ and its conditional expectation on
$\CZ_{k-1}\type{k-1}$ have the same norm in $L^2(\mu)$ and it follows that
they  are equal and that $A$ is measurable with respect to
$\CZ_{k-1}\type{k-1}$.

The announced minimality property of $\CZ_{k-1}$ can be proven by a
method similar to the one used in the proof of~\ref{it:Zk3}.
\end{proof}

\begin{corollary}
\label{cor:Zk}
$Z_0$ is the trivial factor of $X$ and $Z_1$ is its  Kronecker factor.
The sequence of factors $\{Z_k\colon k\geq 0\}$ is increasing.
\end{corollary}
We therefore have a chain of factor maps:
\begin{equation}
\label{eq:chain}
X\to\dots \to Z_{k+1}\to Z_k\to Z_{k-1}\to\dots\to Z_1\to Z_0\ .
\end{equation}

\begin{proof}
All these properties follow from part~\ref{it:Zk2} of the Proposition
by using  the formula $\nnorm f_1=|\int f\,d\mu|$ for the first one,
Lemma~\ref{lemma:mu2} for the second one and the ordering of the
seminorms (Subsection~\ref{subsec:constsemi}) for the last one.
\end{proof}

\subsection{Systems of order $k$}
\begin{definition}
Let $k\geq 0$ be an integer.
A \emph{system of order $k$} is an ergodic system $X$ with $Z_k(X)=X$.
\end{definition}
It is easy to check that for any ergodic system $X$ we have
$Z_k(Z_k(X))=Z_k(X)$ and thus that $Z_k(X)$ is a system of order
$k$.  By Corollary~\ref{cor:Zk} there exists a unique system of
order  zero, the trivial system. The systems of order $1$  are those
which are equal to their Kronecker factor, that is the ergodic
rotations (see Subsection~\ref{subsec:kro}). Every system of order
$k$ is also a system of order $\ell$ for every $\ell >k$.

\subsubsection{Reduction to systems of order $k$.}
\label{subsec:reduclevelk}
We explain here how it suffices to prove the convergence theorems for
systems of order $k$ (for some $k$).
Let $(X,\mu,T)$ be an ergodic system, $f_1,f_2$, $\dots,f_k$ be bounded
functions on $X$ and consider the averages
\begin{equation*}
\tag{\ref{eq:AP}}
\frac 1N\sum_{n=0}^{N-1}f_1(T^nx)f_2(T^{2n}x)\dots f_k(T^{kn}x)
\end{equation*}
as in Theorem~\ref{th:convprog}. Fix $i\in\{1,2,\dots,k\}$. By
part~\ref{it:Zk2} of Proposition~\ref{prop:propertiesZk},
$\nnorm{f_i-\E(f_i\mid\CZ_{k-1})}_k=0$ and thus by
Proposition~\ref{prop:semiprog} the difference between the
averages~\eqref{eq:AP} and the same averages with $\E(f_i\mid\CZ_{k-1})$
substituted for $f_i$ converge to zero in $L^2(\mu)$. In
order to prove Theorem~\ref{th:convprog}, we can thus  assume without
loss that all the
functions $f_i$ are measurable with respect to $\CZ_{k-1}$. Therefore
we can assume that the functions are defined on the associated factor
system $Z_{k-1}(X)$:  we say that this factor is a
\emph{characteristic factor} for the convergence of the
averages~\eqref{eq:AP}. As $Z_{k-1}(X)$ is a system of order $k-1$, it
is sufficient to prove the convergence of these averages under the
additional hypothesis that $X$ is a system of order $k-1$.

The same method applies to the polynomial averages of
Theorem~\ref{th:pol}: for every polynomial family $p_1,p_2,\dots,p_k$
there exists by Proposition~\ref{prop:semipol} an integer $\ell$
such that it suffices to show the convergence under the additional
hypothesis that the system is of order $\ell$.

\subsubsection{Reduction to nilsystems}
\label{subsec:reducnil}
In the rest of the paper we establish a  relation between the
systems of order $k$ and the $k$-step nilsystem described in
Section~\ref{sec:nil}. We need a definition.
\begin{definition}
For each integer $i\geq 1$ let $(X_i,\mu_i,T)$ be a factor of the
system $(X,\mu,T)$ and  assume that this sequence is increasing,
meaning that the sequence $\{ \CX_i \}$ of associated sub-$\sigma$-algebras of
$\CX$ is increasing.
We say that $X$ is the \emph{inverse limit} (or the projective limit)
of the sequence $\{ X_i \}$ if $\CX=\bigvee_i\CX_i$, that is, if
$\CX=\bigcup_i \CX_i$ up to null sets.
\end{definition}

\begin{ST}[\cite{HK4}, Theorem~10.1]
For every $k\geq 1$, every system of order $k$ is an inverse limit of
a sequence of $k$-step nilsystems.
\end{ST}

The convergence results (Theorem~\ref{th:convprog} and~\ref{th:pol})
follow easily since they hold for nilsystems
(Corollaries~\ref{cor:convprognil} and~\ref{cor:convpolnil}) and pass
to inverse limits.

The proof of the Structure Theorem is the longest and the most technical
part of the proofs of convergence and in the next section we can only give a
relatively vague idea of the strategy .
To understand why this proof is
long it is perhaps interesting  to compare the two notions involved
in the theorem. Systems of order $k$ are defined in terms of abstract
ergodic theory, without any mention of a topology or a differentiable
structure on the space. The unique ingredients of this \emph{poor structure}
are a probability space and a measure preserving transformation. On the
other hand, nilsystems are defined in terms of Lie Groups and clearly have
  a rich structure.  Proving the Structure Theorem therefore forces
  us to build this rich structure from scratch.

\subsection{Notes on Section~\ref{sec:factors}}

\subsubsection{Characteristic factors}
The notion of a characteristic factor was already used (without the name) in
Furstenberg's original paper~\cite{F2} and applies to  a wide range of
questions.
Assume for example that we are dealing with the limit behavior of a sequence
of averages depending on some functions.  We say that a factor $Y$ of $X$
is characteristic if the difference between the given averages and the
averages with each function replaced by its conditional expectation on
$\CY$ converges to zero (for the notion of convergence in question).
We are therefore left with studying the limit behavior with $\CY$
substituted for the given system and this can be much easier if $\CY$
has a ``rich'' structure.

Clearly characteristic factors are not unique: a factor containing a
characteristic one is characteristic too. It can be proven that the
factors $Z_k$ defined here are the smallest possible for our
convergence problems. Furstenberg used the much larger
\emph{maximal distal factor}. The ``structure'' of this factor is much weaker
than that of a nilmanifold but is sufficiently rich to make the proof of
Furstenberg's Theorem possible.

\subsubsection{The case of a nilsystem}

The factors $Z_k(X)$ can be described explicitly for nilsystems.

Let  $(X,\mu,T)$ be an ergodic $\ell$-step nilsystem.
We use the notation of  Section~\ref{sec:nil}
and assume moreover that the group $G$ is spanned by its connected component
of the identity and the element $t$ defining the transformation $T$
(it is always possible to reduce to this case).
For every $k\geq 1$,
$(\Lambda  G_{k+1})/G_{k+1}$ is a discrete and cocompact subgroup of
the nilpotent Lie group $G/G_{k+1}$, and $Z_k(X)$ is the nilsystem
$$
Z_k(X)=\frac{G/G_{k+1}}{(\Lambda  G_{k+1})/G_{k+1}}= \frac G{G_{k+1}\Lambda}
$$
endowed with  translation by the projection of $t$ on $G/G_{k+1}$.
This result was already proved by Parry~(\cite{Pa1}) and
Leibman~(\cite{Lei2}) for the Kronecker factor corresponding to the
case $k=1$. Using this
formula with $k=\ell$ we get that every ergodic $k$-step nilsystem is a system
of order $k$.

\section{On the way to the Structure Theorem}
\subsection{A group associated to an ergodic system}
To every ergodic system $(X,\mu,T)$ we associate a group $\CG(X)$ of
measure preserving transformations. The strategy consists in showing that for
sufficiently many systems of order $k$ this group is a nilpotent Lie
group and acts ``transitively'', so that the system can be given the
structure of a nilsystem. We need some notation.

Let $g$ be a measure preserving transformation on $X$ written
$x\mapsto g\cdot x$ and let $k\geq 1$ be an
integer. For each side $\alpha$ of $\{0,1\}^k$ we define the
transformation $g\type k_\alpha$ of $X\type k$ by
$$
(g\type k_\alpha\cdot\bx)_\eps=
\begin{cases}
g\cdot x_\eps&\text{if }\eps\in\alpha\ ;\\
x_\eps & \text{otherwise.}
\end{cases}
$$
This definition coincides with that of the side transformation in the
case that $g=T$.
\begin{definition}
$\CG(X)$ is the group of measure preserving transformations of $X$ such
that for every integer $k\geq 1$ and every side $\alpha$ of
$\{0,1\}^k$, the transformation $g\type k_\alpha$  leaves the measure
$\mu\type k$ invariant.
\end{definition}
The proofs of the following results can be found in Section~5
of~\cite{HK4}.
We write $\CG$ instead of $\CG(X)$ except when some ambiguity can occur.
This group is a Polish group when endowed with the topology of convergence
in probability. It contains $T$ and every measure preserving
transformation of $X$ commuting with $T$.

\begin{lemma}
\label{lemma:inducedG}
Let $g\in\CG(X)$. Then for every $k$ the transformation $g$ of $X$  maps the
factor $\CZ_k$ to itself and thus induces a transformation of $Z_k$,
which belongs to $G(Z_k)$.
\end{lemma}
\begin{lemma}
If $X$ is a system of order $k$ then $\CG(X)$ is a $k$-step nilpotent
group.
\end{lemma}
\begin{proof}
Let $g_1,g_2,\dots,g_{k+1}\in \CG$ and
$h=[[\dots[g_1,g_2],g_3],\dots,g_{k+1}]$.

Let $\alpha_1,\alpha_2,\dots\alpha_{k+1}$ be the $k+1$ sides of $\{0,1\}^k$
containing $\zero$. The measure $\mu\type {k+1}$ is invariant under each
transformation $g_{\alpha_i}\type {k+1}$ and thus also under their commutator
$[[\dots[g_{\alpha_1}\type {k+1},g_{\alpha_2}\type {k+1}],
g_{\alpha_3}\type {k+1}],\dots,g_{\alpha_{k+1}}\type {k+1}]$. But is 
is easy to check that this
transformation is equal to the transformation $h\type {k+1}_\zero$,
given by $(h\type {k+1}_\zero\cdot \bx)_\zero=h\cdot x_\zero$ and
$(h\type {k+1}_\zero\cdot \bx)_\eps=x_\eps$ for
$\eps\neq\zero$.

Let $A$ be a subset of $X$. As $X$ is of level $k$, $Z_k(X)=X$ and by
definition there exists a subset $B$ of $X\typ{k+1}$ with
$\one_A(x_\zero)=\one_B(\tilde\bx)$ for $\mu\type{k+1}$ almost
every $\bx$. Applying the transformation $h\type {k+1}_\zero$ we get
that $\one_A(h\cdot x_\zero)=\one_B(\tilde\bx)$ ($\mu\type{k+1}$
a.e.) and thus that
$\one_A(x_\zero)=\one_A(h\cdot x_\zero)$ ($\mu$-a.e.). This shows
that every subset $A$ of $X$ is invariant under $h$ and that $h=\id$.
\end{proof}

\subsection{Relations between two consecutive factors}

\begin{definition}
Let $(Y,\nu,T)$ be a system, $U$ a compact abelian group endowed with
its Haar measure $m_U$ and $\rho\colon Y\to U$ a measurable map. Let
$X=Y\times U$ be endowed with the measure $\mu=\nu\times m_U$ and
with the transformation $T$ given by $T(y,u)=(Ty,u+\rho(y))$. Then we
say that $X$ is an \emph{extension of $Y$ by $U$} and $\rho$ is
called the \emph{cocycle defining the extension.}
\end{definition}
We note that $Y$ is a factor of $X$, with the factor map given by
$(y,u)\mapsto y$. For each $v\in U$, the transformation
$R_u\colon (y,u)\mapsto(y,u+v)$ of $X$ preserves $\mu$ and commutes
with $T$; it is called a \emph{vertical rotation}.
\begin{proposition}
\label{prop:ext}
Let $(X,\mu,T)$ be a system of level $k$ and
$(Y,\nu,T)=Z_{k-1}(X)$. Then
$X$ is an extension of $Y$ by a compact abelian group $U$.
Moreover, the cocycle $\rho\colon Y\to U$ defining this
extension satisfies:

There exists a map $F\colon Y\type k\to U$ with
\begin{equation}
\label{eq:typek}
\sum_{\eps\in\{0,1\}^k}
(-1)^{\lvert\eps\rvert}
\rho(y_\eps)=F(T\type k\by)-F(\by)
\end{equation}
for $\nu\type k$-almost every $\by\in Y\type k$.
\end{proposition}
A cocycle $\rho$ satisfying the functional equation~\eqref{eq:typek} for
some map $F$ is called \emph{a cocycle of type $k$}.
\begin{proof}[Idea of the proof]
The first step of the proof
uses the notion of  an isometric extension as defined by
Furstenberg (see~\cite{F3}). Part~\ref{it:Zk1} of
Proposition~\ref{prop:propertiesZk}  implies that the
$T\type k$-invariant $\sigma$-algebra $\CI\type k$ of
$X\type k$ is measurable with respect to $\CW\type k$, where $W$ is
some factor of $X$ which is an isometric extension of $Y$.
The minimality property~\ref{it:Zk4} of the same Proposition then
gives that $X=W$, that is, $X$ is an isometric extension of $Y$.

In particular there exists a compact group $U$, acting on $X$ by
measure preserving transformations and inducing the trivial
transformation on $Y$. It is then proven that this group of
transformations is included in the center of $\CG(X)$ and in particular
is an abelian group; this shows that $X$ is an extension of $Y$
by the compact abelian group $U$. We identify $X$ with $Y\times
U$. Let $\rho\colon Y\to U$ be the cocycle defining
this extension.

Let $\chi$ be a character of $U$, that is a continuous group
homomorphism from $U$ to the circle group $\CS^1$. Let $\phi$ be
the function on $X=Y\times U$ defined by $\phi(y,u)=\chi(u)$
and $\Phi$ the function defined on $X\type k=Y\type k\times
U\type k$ by
$$
\Phi(\by,\bu)=
\prod_{\eps\in\{0,1\}^k}
C^{\lvert\eps\rvert} \phi(x_\eps)=
\chi\Bigl(\sum_{\eps\in\{0,1\}^k}
(-1)^{\lvert\eps\rvert}
u_\eps\Bigr)
$$
for $\bx=(\by,\bu)$ with $\by\in Y\type k$ and $\bu\in U\type k$.

As $X$ is of type $k$, $\nnorm\psi_{k+1}\neq 0$ by part~\ref{it:Zk2}
of Proposition~\ref{prop:propertiesZk}. By construction of the
seminorms this means that the function $\Psi=\E(\Phi\mid\CI\type k)$
is not identically zero. This function is invariant under $T\type k$
and satisfies $\Psi(R_u\bx)=\chi(u)\Psi(\bx)$ for every $u\in U$
and $\mu$-almost every $x\in X$. By using the ergodicity property
(Lemma~\ref{lemma:ergod}) of $\mu\type  k$ it can be showed that
there exists a function with the same properties and everywhere
nonzero, and thus a function $F_\chi$ of modulus $1$ with the same
properties. The invariance of this function gives:
$$
\chi\Bigl(\sum_{\eps\in\{0,1\}^k}
(-1)^{\lvert\eps\rvert}
\rho_k(y_\eps)\Bigr)=F_\chi(T\type k\by)\,F_\chi(\by)\inv\ .
$$
The existence of a function $F_\chi$ with this property for every
character $\chi$ of $U$ implies the existence of a function $F$
satisfying ~\eqref{eq:typek} by classical results
about cocycles.
\end{proof}

\subsection{A technical tool}
 From this point, the proof of the Structure Theorem proceeds by induction
on $k$. We give here a technical tool used in the induction.

For every $s$ in the Kronecker factor $Z_1$ we have defined in
Subsection~\ref{subsec:mu2} a
measure $\mu_s$ on $X\times X$, invariant under $T\times T$. For
almost every $s$, the system $(X\times X,\mu_2,T\times T)$ is ergodic
and we denote it by $X_s$.
\begin{proposition}
Let the hypotheses and the notation be as in
Proposition~\ref{prop:ext}.
Then, for almost every $s\in Z_1$,
$X_s$ is a system of level $k$,
$Y_s$ a system of level $k-1$ and
$X_s$ is an extension of $Y_s$ by the compact abelian group
$U\times U$, given by the
cocycle $(y_0,y_1)\mapsto (\rho(y_0),\rho(y_1))$ which is of type $k$.

Moreover, $Z_{k-1}(X_s)$ is an extension of $Y_s$ by $U$,
given by the
cocycle $(y_0,y_1)\mapsto \rho(y_0)-\rho(y_1)$ which is of type $k-1$.
\end{proposition}
This proposition  plays the role of a
stepladder, allowing us to climb from one level to the next one:
assume that some properties have been shown for systems of order
$k-1$ and let $X$ be a system of order $k$. Let $U$ and $\rho$ be as
in Proposition~\ref{prop:ext}. Then we can use these properties for
the systems $Z_{k-1}(X_s)$ and this gives  information on the group
$U$ and the cocycle $\rho$. This method is used in particular to prove:
\begin{lemma}
Let $X,Y,U$ and $\rho$ be as in Proposition~\ref{prop:ext}. Then $U$
is connected.
\end{lemma}

\subsection{Toral systems}
Recall that a compact abelian group can be given a structure of Lie
group if and only if the connected component of the  its unit element
is a finite dimensional torus. This property is equivalent to saying that
its dual group is finitely generated. It follows that every compact
(metrizable) abelian group can be represented as an inverse limit of
a sequence of compact abelian Lie groups. In particular every compact
connected abelian group can be represented as an inverse limit of a
sequence of finite dimensional tori. This motivates the next
definition.
\begin{definition}
Let $k\geq 1$ be an integer and $(X,\mu,T)$ be a system of order $k$.
We say that this system is \emph{toral} if its Kronecker factor
$Z_1(X)$ is a compact abelian Lie group and if $Z_j(X)$ is an extension
of $Z_{j-1}(X)$ by a finite dimensional torus for $2\leq j\leq k$.
\end{definition}
The structure Theorem can be split into the next two Propositions:
\begin{proposition}
For every integer $k\geq1$, every system of order $k$  is
the inverse limit of a sequence of toral systems of
order $k$.
\end{proposition}
\begin{proposition}
Every toral system of order $k$ is  $k$-step nilsystem.
\end{proposition}
More precisely, $\CG(X)$ is a Lie group and acts transitively on $X$.
$X$ can be identified with $G(X)/\Lambda$ where $G(X)$ is the
subgroup of $\CG(X)$ spanned by the connected component of its
identity and $T$. We give an idea of the method of the proof.

Let $X$ be a system of order $k$; we use the notation of
Proposition~\ref{prop:ext}.  As we assume that the result holds for
systems of order $k-1$ it holds in particular for $Y=Z_{k-1}(X)$. Each
element $h$ of $G(Y)$ is lifted to an element $\tilde h$ of $G(X)$; 
this
means that $h$ is the transformation induced by $\tilde h$ on $Y$ as
in Lemma~\ref{lemma:inducedG}. This is done first for $h$ belonging
to $G(Y)_{k-1}$ then for $h\in G(Y)_{k-1}$, and so on, following the
lower central series of $G(Y)$ upwards. Each step uses the functional
equation~\eqref{eq:typek} satisfied by the cocycle $\rho$.

\appendix
\section*{Appendix: Further comments }
\setcounter{section}{1}
\subsection{The cubic averages}
\label{note:cubic}
The paper~\cite{HK4} also contains  the proof of convergence of another
type of multiple ergodic average, the \emph{cubic averages}, already
proven by Bergelson for the cubes of dimension $2$:
\begin{theorem*}[Bergelson~\cite{Berg2}]
Let $(X,\mu,T)$ be a system and let $f,g,h$ be three bounded
functions on $X$. Then the averages
$$
   \frac 1{N^2} \sum_{m,n=0}^{N-1}f(T^mx)g(T^nx)h(T^{m+n}x)
$$
converge in $L^2(\mu)$ as $N\to+\infty$.
\end{theorem*}
For the general case we use here the notation introduced at the top of
Section~\ref{sec:seminorms}.

\begin{theorem}[Host \& Kra~\cite{HK4}]
\label{th:cubes}
Let $(X,\mu,T)$ be a system and
let $f_\eps$, $\eps\in\{0,1\}^k\setminus\{\zero\}$ be $2^k-1$
bounded functions on $X$. Then the averages
\begin{equation}
\label{eq:avcubes}
\frac 1{N^k}\sum_{n_1,n_2,\dots,n_k=0}^{N-1}\;
\prod_{\substack{\eps\in\{0,1\}^k\\ \eps\neq\zero}}
f_\eps(T^{\bn\cdot\eps}x)
\end{equation}
converge in $L^2(\mu)$.
\end{theorem}

The strategy for the proof of this theorem is
the same as for the averages along arithmetic progressions and for the
polynomial averages.
Clearly it suffices to prove the result for ergodic systems. First, a 
bound similar to~\eqref{eq:semiprog} holds for
the cubic averages, with a similar proof using a
multidimensional version of the van der Corput Lemma. Here $(X,\mu,T)$
is an ergodic system.

\begin{proposition*}
Let $f_\eps$, $\eps\in\{0,1\}^k\setminus\{\zero\}$ be $2^k-1$
functions on $X$ with $\nnorm{f_\eps}_\infty\leq 1$ for every
$\eps$.
   The $\limsup$ of the norm in $L^2(\mu)$ of the
averages~\eqref{eq:avcubes} is bounded by
$\displaystyle\min_\eps\nnorm{f_\eps}_k$.
\end{proposition*}

By the same arguments as in Subsections~\ref{subsec:reduclevelk}
and~\ref{subsec:reducnil} it is therefore possible to restrict to the
case that $X$ is a  $(k-1)$-step nilsystem. Then
Theorem~\ref{th:cubes} follows
from a generalization of Theorems~\ref{th:nilsystems1} 
and~~\ref{th:nilsystems2}
and of Corollary~\ref{cor:nilsystems1} to the case of several
commuting translations on a nilmanifold.

The cubic averages are directly linked to the measures $\mu\type k$:
\begin{proposition*}
Let $f_\eps$, $\eps\in\{0,1\}^k$, be $2^k$ functions  on
$X$. Then
$$
\frac 1{N^k} \sum_{n_1,n_2,\dots,n_k=0}^{N-1}\;
\int \prod_{\eps\in\{0,1\}^k}
f_\eps(T^{\bn\cdot\eps}x)\,d\mu(x)
\to
\int\prod_{\eps\in\{0,1\}^k}f_\eps(x_\eps)\,
d\mu\type k(\bx) \ .
$$
\end{proposition*}
By using this result with all functions equal to $\one_A$ and the
inequality $\nnorm{\one_A}_k\geq\nnorm{\one_A}_1= \mu(A)$ we get:
\begin{corollary*}
For every subset $A$ of $X$,
$$
\lim_{N\to+\infty}\frac 1{N^k} \sum_{n_1,n_2,\dots,n_k=0}^{N-1}\;
\mu\Bigl(\bigcap_{\eps\in\{0,1\}^k}T^{\bn\cdot\eps}A
\Bigr)\geq\mu(A)^{2^k}\ .
$$
\end{corollary*}

The same results hold when the averages on $[1,N)^k$ are replaced by
averages on a sequence of parallelepipeds whose minimal size tend to
infinity, or more generally by averages on a F\o lner sequence. This
version of the last Corollary has a combinatorial interpretation in
terms of sets of integer of positive upper density  but we do not state it
here (see~\cite{HK4}).

The convergence a.e. of the averages~\eqref{eq:avcubes} has been
recently proven by Assani (\cite{A}).

\subsection{Gowers norms and their relations to arithmetic
progressions}
\label{subsec:gowersnorms}
\subsubsection{The definition}
When writing our paper~\cite{HK4} we became aware of Gowers'
paper~\cite{Gowers1} where he introduced some norms very similar to
our seminorms. It turns out that they are  identical to our seminorms
when computed in the particular case that $X=\Z/N\Z$ endowed with the
uniform measure and
with the transformation $x\mapsto x+1\bmod N$. These norms were extensively
used by Green and Tao~(\cite{GT},~\cite{Tao}).
We recall here their definition, with notation modified in order to
fit  with ours.

Let $N\geq 2$ be an integer and let $G=\Z/N\Z$ be endowed with its
normalized Haar measure $m$. Let $\CC(G)$ denotes the space of complex valued
functions on $G$.
For $f\in \CC(G)$,  $\nnorm f_1$  is defined to be
$\bigl|\int f(x)\,dm(x)\bigr|$.
For $t\in G$ let $f_t$ be the function $x\mapsto f(x+t)$ and
define by induction
\begin{equation}
\label{eq:defsemianal}
\text{for }k\geq 1,\
\nnorm f_{k+1}
=\Bigl(\int\nnorm{f.\overline{f_t}}_k^{2^k}\,dt\Bigr)^{1/2^{k+1}}\ .
\end{equation}
$\nnorm f_k$ can also be defined by a closed formula. 
For $\bt=(t_1,t_2,\dots,t_k)\in G^k$ and $\eps\in\{0,1\}^k$ we write
$\eps\cdot\bt=\eps_1t_1+\eps_2t_2+\dots+\eps_kt_k$.
For $f\in\CC(G)$ and $k\geq 1$ we have by induction
\begin{equation}
\label{eq:defsemianal2}
\nnorm f_k
=\Bigl(\int\prod_{\eps\in\{0,1\}^k}C^{\lvert\eps\rvert}f(x+\eps\cdot\bt)\,
dm(x)\,dm(t_1)\,dm(t_2)\,\dots\,dm(t_k)\Bigr)^{1/2^k}\ .
\end{equation}
An inequality similar to the Cauchy-Schwarz Inequality follows:
when $f_\eps$, $\eps\in\{0,1\}^k$, are functions on $G$,
\begin{equation}
\label{eq:cauchyanal}
\Bigl|
\int\prod_{\eps\in\{0,1\}^k}
f_\eps(x+\eps\cdot\bt)\,dm(x)\,dm(t_1)\,dm(t_2)\,\dots\,dm(t_k)
\Bigr|
\leq \prod_{\eps\in\{0,1\}^k}\nnorm {f_\eps}_k\ .
\end{equation}
This implies that each map $f\mapsto\nnorm f_k$ is subadditive and hence
a seminorm. For $k=2$, relation~\eqref{eq:defsemianal2} gives:
$$
\nnorm f_2= \Bigl(\sum_{j\in G}\bigl|\widehat f(j)\bigr|^4\Bigr)^{1/4}
$$
where $\widehat f$ is the Fourier Transform of $f$, defined on the
dual group $\widehat G=G$ of $G$. This expression can only be zero
when $f$ is the zero function and thus  $\nnorm\cdot_2$ is a
norm on $\CC$. By induction and by using
definition~\eqref{eq:defsemianal}, it can be checked that for every
$f\in\CC(G)$ we
have $\nnorm f_{k+1}\geq\nnorm f_k$ for every $k$ and thus
$\nnorm \cdot_k$ is a norm on $\CC(G)$ for every $k\geq 2$.

\subsubsection{Using Gowers norms}
These norms are used by Gowers, Green and Tao to control the
integral~\eqref{eq:SzAnal}
appearing in the analytic form of Szemer\'edi's Theorem.
\begin{proposition}
\label{prop:Szanalprog}
Let $N\geq 2$ be an integer and let $\Z/N\Z$ be endowed with its
normalized Haar measure. Let $\ell\geq 2$ be an integer and
$f_0,f_1,f_2,\dots,f_{\ell-1}$ be functions on  $\Z/N\Z$ with $|f_i|\leq 1$
for $0\leq i\leq\ell-1$. Then
\begin{multline}
\label{eq:Szanalprog}
\Bigl| \iint
f_0(x)f_1(x+y)f_2(x+2y)\dots f_{\ell-1}(x+(\ell-1)y)\,dm(x)\,dm(y)
\Bigr|\\
\leq \min_{0\leq i\leq\ell-1}\nnorm{f_i}_{\ell-1}\ .
\end{multline}
\end{proposition}
The starting point of Gowers' method (\cite{Gowers1}) can be summarized
very roughly as follows.
Let $A$ be a subset of $G$, $f=\one_A$ and $g=f-m(A)$.
He distinguishes two cases: If $\nnorm g_{\ell-1}$ is small then the
integral in~\eqref{eq:SzAnal} is close to the same integral
with the constant $m(A)$ substituted for $f$ and thus is large. If  $\nnorm
g_{\ell-1}$ is large then he shows that the restriction of $f$
to some relatively large subset of $G$ has some strong arithmetic
properties and behaves in some respects like a polynomial; he deduces that
the integral~\eqref{eq:SzAnal} is large in this case also.

The way  Tao (\cite{Tao})  and Green and Tao (\cite{GT}) use Gowers' norms
is much closer to  ergodic theory and in particular to the way we use
the seminorms here.  We try here to make this
similarity apparent but the reader must be
advised that the contents of the next few lines is nothing more than 
an oversimplification.

Tao decomposes  any function $f$
on $G$ as a sum of a function $g$ with small norm and a ``anti-uniform''
function $h$. The contribution of $g$ to the
integral in~\eqref{eq:SzAnal} is small and the function $h$ can be
viewed as the conditional expectation of $f$ relatively to some
$\sigma$-algebra, comparable to the factor used in~\cite{FKO} except
that it is not invariant under translation.
The contribution of $h$
to the integral is bounded from below by using van der Waerden's Theorem.

The theorem of Green and Tao (\cite{GT}) about the existence of
arithmetic progressions in primes seems completely out of the range of
ergodic theory because the primes have zero density in the integers and
thus the Correspondence Principle does not apply. However the strategy
is similar.
  They show that in the preceding decomposition the uniform
norm of the function $h$ can be bounded independently of any uniform
bound of the function $f$, assuming only that this function is bounded
by a ``quasirandom'' function. Therefore the analytic form of
Szemer\'edi's Theorem can be used to show that the function $h$ gives
a ``large'' contribution to the integral in~\eqref{eq:SzAnal}.
Moreover the function $g$ also is bounded by a quasirandom function
and its contribution to the integral is small, due to an extension of
Proposition~\ref{prop:Szanalprog} to this case. They therefore get
a generalization of the analytic form of Szemer\'edi's Theorem under
the weaker hypothesis. This result is then used for
a function closely related to the indicator function of primes.

The  decomposition used in both cases is parallel to the decomposition
$$
f=\E(f\mid\CZ_{k-1})+\bigl(f-\E(f\mid\CZ_{k-1})\bigr)
$$
used in subsection~\ref{subsec:reduclevelk} of this paper but the
authors do not need a precise description of the ``factor'' comparable
to the Structure Theorem and do not use the machinery of nilpotent
groups.  It can be conjectured that there exists a hidden link between
the combinatorial constructions and the nilpotent groups and we
believe that making this link explicit is a very interesting challenge.


\begin{thebibliography}{FKO82}

\bibitem[As]{A}
I.~Assani.
Pointwise convergence of averages along cubes.
{\em Preprint} (2004).


\bibitem[AuGH]{AGH}þ
L. Auslander, L. Green and F. Hahn.
Flows on homogeneous spaces.
{\em Ann. Math. Studies}
{\bf 53}, Princeton Univ. Press (1963).

\bibitem[B1]{Berg1}
V.~Bergelson.
Weakly mixing PET.
{\em Erg. Th. \& Dyn. Sys.}, {\bf 7}  (1987), 337--349.

\bibitem[B2]{Berg2}
V.~Bergelson.
The multifarious Poincar\'{e} Recurrence Theorem.
{\em Descriptive Set Theory and Dynamical Systems}, Eds. M. Foreman,
A.S. Kechris, A. Louveau, B. Weiss.  Cambridge University Press, New York
(2000), 31-57.

\bibitem[BHK]{BHK}
V.~Bergelson, B.~Host and B.~Kra, with an Appendix by I. Ruzsa. 
Multiple recurrence and nilsequences. {\em Inventiones Math.}
{\bf 160} (2005), 261--303. 

\bibitem[BL96]{BL}
V.~Bergelson and A.~Leibman.
Polynomial extensions of van der Waerden's and Szemer\'edi's
Theorems.
{\em Journal Amer. Math. Soc.}, {\bf 9} (1996), 725--753.

\bibitem[BMc]{BMC}
V.~Bergelson and R.~McCutcheon.
An Ergodic IP Polynomial Szemer\'edi Theorem.
{\em Mem. Amer. Math. Soc.} {\bf 146} (2000), \#146.

\bibitem[Bo1]{Bourgain1}
J.~Bourgain.
Pointwise ergodic theorems for arithmetic sets.
{\em Inst. Hautes \'{E}tudes Sci. Publ. Math.}, {\bf 69}
    (1989), 5--45.

\bibitem[Bo2]{Bourgain2}
J.~Bourgain.
On triples in arithmetic progression. {\em Geom. Funct. Anal.} {\bf
9}, 1999, 968--984.

\bibitem[CL1]{CL1}
J.-P.~Conze and E.~Lesigne.
Th\'{e}or\`{e}mes ergodiques pour des mesures diagonales.
{\em Bull. Soc. Math. France}, \textbf{112} (1984), 143--175.

\bibitem[CL2]{CL3}
J.-P. Conze and E.~Lesigne.
Sur un th\'{e}or\`{e}me ergodique pour des mesures diagonales.
{\em Publications de l'Institut de Recherche de Math\'ematiques de
    Rennes, Probabilit\'es}, 1987.

\bibitem[CL3]{CL2}
J.-P. Conze and E.~Lesigne.
Sur un th\'{e}or\`{e}me ergodique pour des mesures diagonales.
{\em C. R. Acad. Sci. Paris, S\'{e}rie {I}},
\textbf{306} (1988), 491--493.

\bibitem[ET]{ET}
P.~Erd\"os \& T.~Tur\'an.
On some sequences of integers.
{\em J. London Math. Soc.} {\bf 11} (1936), 261--264.

\bibitem[FK1]{FK1}
N.~Frantzikinakis and B.~Kra.
Polynomial averages converge to the product of integrals.
 {\em Israel J. of Maths.} {\bf 148} (2005), 267--276.

\bibitem[FK2]{FK2}
N.~Frantzikinakis and B.~Kra.
 Convergence of multiple ergodic averages for some commuting 
ransformations. 
{\em Erg. Th. \& Dyn. Sys.} {\bf 25} (2005), 799--809. 


\bibitem[Fr1]{Fr1}
G.~A.~Freiman.
Foundations of a Structural Theory of Set Addition.
{\em Translations of Mathematical Monographs} {\bf 37} (1973), Amer.
Math. Soc.,
Providence.

\bibitem[Fr2]{Fr2}
G.~A.~Freiman.
Structure theory of set addition.
{\em Structure theory of set addition}, Asterisque {\bf 258} (1999), 1--33.

\bibitem[Fu1]{F1}
H.~Furstenberg.
Strict ergodicity and transformations of the torus.
{\em Amer. J. of Mathematics}, {\bf 83} (1961), 573--601.

\bibitem[Fu2]{F2}
H.~Furstenberg.
Ergodic behavior of diagonal measures and a theorem of
Szemer\'edi  on arithmetic progressions.
{\em J. d'Analyse Math.}, \textbf{31} (1977), 204--256.

\bibitem[Fu3]{F3}
H.~Furstenberg.
{\em Recurrence in Ergodic Theory and Combinatorial Number
Theory.}
Princeton Univ. Press (1981).

\bibitem[FuK1]{FK79}
H.~Furstenberg and Y.~Katznelson.
An ergodic Szemer\'edi theorem for commuting transformations.
{\em J. Analyse Math.}  {\bf 34} (1978), 275--291.

\bibitem[FuK2]{FK91}
H.~Furstenberg and Y.~Katznelson.
A density version of the Hales-Jewett theorem.
{\em J. Anal. Math.}  {\bf 57}  (1991), 64--119.

\bibitem[FuKO]{FKO}
H.~Furstenberg, Y.~Katznelson and D.~Ornstein.
The ergodic theoritical proof of Szemer\'edi's theorem.
{\em Bull. Amer. Math. Soc.} {\bf 7} (1982), 527--552.

\bibitem[FuW]{FW}
H.~Furstenberg and B.~Weiss.
A mean ergodic theorem for
$\frac{1}{N} \sum_{n=1}^n f({T}^nx)g({T}^{n^2}x)$.
{\em Convergence in Ergodic Theory and Probability},
Eds.: Bergelson, March, Rosenblatt. Walter de Gruyter \& Co, Berlin, New
York (1996), 193--227

\bibitem[G1]{Gowers1}
T.~Gowers.
A new proof of {S}zemer\'edi's theorem.
{\em Geom. Funct. Anal.}, {\bf 11} (2001), 465-588.

\bibitem[G2]{Gowers2}
T.~Gowers.
Hypergraph regularity and the multidimensional Szemer\'edi Theorem.
{\em Preprint} (2004).

\bibitem[GrT]{GT}
B.~Green and T.~Tao.
The primes contain arbitrarily long arithmetic progressions.
{\em Preprint} available at http://arxiv.org (math.NT/0404188) (2004).

\bibitem[HK1]{HK1}
B.~Host and B.~Kra.
Convergence of Conze-Lesigne Averages.
{\em Erg. Th. \& Dyn. Sys.}, \textbf{21} (2001), 493--509.

\bibitem[HK2]{HK4}
B.~Host and B.~Kra.
Nonconventional ergodic
averages and nilmanifolds. 
{\em Annals of Math.} {\bf 161} (2005), 397-Ð488. 

\bibitem[HK3]{HK6}
B.~Host and B.~Kra.
Convergence of polynomial ergodic averages. 
{\em Isr. J. Math.} {\bf 149} (2005), 1--9. 


\bibitem[HK4]{HKcubes}
     B.~Host and B.~Kra.
Averaging along cubes.
{\em Dynamical Systems and Related Topics},
Eds. Brin, Hasselblatt, Pesin. Cambridge University Press, Cambridge
(2004).

\bibitem[Kh]{Khintchine}
A. Y.~Khintchine.
Eine Versch\"{a}rfung des Poincar\'{e}schen
"Wiederkehrsatzes".
{\em Comp. Math.}, \textbf{1} (1934), 177--179.



\bibitem[L2]{Lei2}
A.~Leibman.
Pointwise convergence of ergodic averages for polynomial sequences 
of translations on a nilmanifold. {\em 
Erg. Th. \& Dyn. Sys.} {\bf 25} (2005), 201--213. 

\bibitem[L3]{Lei3}
A.~Leibman.
Convergence of multiple ergodic averages along polynomials of several 
variables. {\em Isr. J. Math.} {\bf 146} (2005), 303Ð-316. 

\bibitem[Le1]{Lesigne1}
E.~Lesigne.
R\'esolution d'une \'equation fonctionnelle.
{\em Bull. Soc. Math. France}, \textbf{112} (1984), 177--19.

\bibitem[Le2]{Lesigne3}
E.~Lesigne.
Th\'eor\`emes ergodiques ponctuels pour des mesures diagonales.
      Cas des syste\`emes distaux.
{\em Ann. Inst. Henri Poincar\'{e}}, {\bf 23} (1987), 593--612.

\bibitem[Le3]{Lesigne4}
E.~Lesigne.
Th\'eor\`emes ergodiques pour une translation sur une
      nilvari\'et\'e.
{\em  Erg. Th. \& Dyn. Sys.}, \textbf{9-1} (1989), 115--126.

\bibitem[Le4]{Le}
E. Lesigne.
Sur une nil-vari\'et\'e, les parties minimales
associ\'ees \`a une translation sont uniquement ergodiques.
{\em Ergod. Th. \& Dyn. Sys.}
{\bf 11} (1991), 379--391.

\bibitem[Le5]{Lesigne2}
E.~Lesigne.
\'Equations fonctionelles, couplages de produits gauches et
      th\'eor\`emes ergodiques pour mesures diagonales.
{\em Bull. Soc. Math. France}, \textbf{121} (1993), 315--351.

\bibitem[M]{Malcev}
A. Malcev.
On a class of homogeneous spaces.
{\em Amer. Math. Soc. Transl.}
{\bf 39} (1951)

\bibitem[P1]{Pa1}
W. Parry.
Ergodic properties of affine transformations and flows on
nilmanifolds.
{\em Amer.  J. Math.}
{\em 91} (1969), 757--771.

\bibitem[P2]{Pa2}
W. Parry.
Dynamical systems on nilmanifolds.
{\em Bull. London Math. Soc.}
{\bf 2} (1970), 37--40.

\bibitem[Ra]{Ratner}
M.~Ratner.
On Raghunathan's measure conjecture.
{\em Ann. Math.}, \textbf{134} (1991), 545-607.

\bibitem[Ro]{Roth}
K.~F.~Roth.
On certain sets of integers.
{\em J. London Math. Soc.} {\bf 28} (1953), 245--252.

\bibitem[Ru]{R}
D.J.~Rudolph.
Eigenfunctions of $T\times S$ and the Conze-Lesigne algebra.
{\em Ergodic Theory and its Connections with Harmonic Analysis},
Eds.: Petersen/Salama, Cambridge University Press, New
York (1995), 369--432.

\bibitem[Sh]{Sh}
N.~Shah.
Invariant measures and orbit closures on homogeneous spaces for
actions of subgroups.
{\em Lie groups and ergodic theory} (Mumbai, 1996)
Tata Inst. Fund. Res., Bombay (1998), 229--271.

\bibitem[Sz1]{Sz1}
E.~Szemer\'edi.
On sets of integers containing no four elements in arithmetic
progression.
{\em Acta Math. Acad. Sci. Hungar.} {\bf 20} (1969), 89--104.

\bibitem[Sz2]{Sz2}
E.~Szemer\'edi.
On sets of integers containing no $k$ elements in arithmetic
progression.
{\em Acta Arith.} {\bf 27} (1975), 199-245.

\bibitem[T]{Tao}
T.~Tao.
A quantitative ergodic proof of Szemer\'edi's theorem.
{\em Preprint} (2004).

\bibitem[Z1]{Ziegler}
T.~Ziegler.
A nonconventional ergodic theorem for a nilsystem.
{\em Erg. The. Dyn. Sys.} {\bf 25} (2005), 1357Ð-1370. 
\end{thebibliography}
\end{document}